\documentclass[preprint,12pt]{article}

\usepackage{amssymb}

\usepackage{pictex}
\usepackage{amsmath}
\usepackage{amscd}
\usepackage{amstext}
\usepackage{amssymb}
\usepackage{amsfonts}
\usepackage{latexsym}
\usepackage[dvips]{color}
\usepackage{tikz-cd}
\newtheorem{T}{Theorem}
\newtheorem{Lem}{Lemma}
\newtheorem{Prop}{Proposition}
\newtheorem{Def}{Definition}
\newtheorem{Cor}{Corollary}
\newtheorem{Rem}{Remark}

\newcommand{\bt}{\begin{T}}
	\newcommand{\bl}{\begin{Lem}}
		\newcommand{\bp}{\begin{Prop}}
			\newcommand{\bc}{\begin{Cor}}
				\newcommand{\bd}{\begin{Def}}
					\newcommand{\br}[2]{\begin{Rem}\label{#1}{\rm #2}}
						\newcommand{\er}{ \end{Rem}}
					\newcommand{\et}{\end{T}}
				\newcommand{\el}{\end{Lem}}
			\newcommand{\ep}{\end{Prop}}
		\newcommand{\ec}{\end{Cor}}
	\newcommand{\ed}{\end{Def}}

\newcommand{\be}{\begin{equation}}
\newcommand{\ee}{\end{equation}}
\newcommand{\beq}{\begin{eqnarray}}
\newcommand{\eeq}{\end{eqnarray}}
\newcommand{\beqq}{\begin{eqnarray*}}
	\newcommand{\eeqq}{\end{eqnarray*}}
\def\supp{\mathop{\rm supp \,}\nolimits}
\def\vol{\mathop{\rm vol}\nolimits}

\def\dist{\mathop{\rm dist \,}\nolimits}

\def\id{\mathop{\rm id \,}\nolimits}
\def\tr{\mathop{\rm tr \,}\nolimits}
\def\ext{\mathop{\rm ext \,}\nolimits}
\def\Id{\mathop{\rm Id}\nolimits}

\newcommand{\reff}[1]{{\rm (\ref{#1})}}
\newcommand{\bpr}{{\bf Proof.}\ }
\newcommand{\epr}{\hspace*{\fill}\rule{3mm}{3mm}\\}

\newcommand{\bspq}{B^s_{p,q}}

\newcommand{\Rd}{{\mathbb R}^{d}}
\newcommand{\R}{{\mathbb R}}
\newcommand{\N}{{\mathbb N}}
\newcommand{\Z}{\mathbb Z}

\newcommand{\C}{{\mathbb C}}

\newcommand{\cl}{{\mathcal{L}}}

\newcommand{\psigma}{\ensuremath{\sigma_{\mathsf p}}}
\newcommand{\esigma}{\ensuremath{\sigma_{\mathsf e}}}
\newcommand{\dom}{\ensuremath{\mathrm{dom}}}
\newcommand{\psigmag}{\ensuremath{\sigma_{\mathsf p}^{\mathsf \gamma}}}

\begin{document}




\title{Some properties of  block-radial functions and Schr\"odinger type operators with block-radial potentials}
\author{Alicja Dota\thanks{Alicja Dota was supported by the Ministry of Science and Higher Education of Poland, Grant No. 04/43/DSPB/0095.}, Leszek Skrzypczak\thanks{Leszek Skrzypczak was supported by National Science Center, Poland, Grant No. 2014/15/B/ST1/00164.} }
\date{ }
\maketitle

%
\begin{abstract}
 Let $ R_\gamma B^{s}_{p,q}(\Rd)$ be a subspace of the Besov space $B^{s}_{p,q}(\Rd)$ that consists of block-radial functions. We prove that   the asymptotic behaviour of the entropy numbers of compact embeddings $\id: \: R_\gamma B^{s_1}_{p_1,q_1}(\R^d) \rightarrow  R_\gamma B^{s_2}_{p_2,q_2}(\R^d)$ depends  on the number of  blocks of the lowest  dimension, the parameters $p_1$ and $p_2$, but is independent of  the smoothness parameters $s_1$, $s_2$. We apply the asymptotic behaviour to estimation of  powers of a negative spectra of Schr\"odinger type operators with block-radial potentials. This part essentially relies on the Birman-Schwinger principle.  
 
\end{abstract}








\section{Introduction}

In recent years, some attention has been paid to describing compactness of embeddings of
function spaces of Besov and Sobolev type by  different quantities, in particular, by corresponding sequences of entropy  and approximation numbers. The study was motivated by the program    formulated  by { D.~Edmunds} and { H.~Triebel}. In \cite{ET} they  proposed  investigation of    spectral properties of certain pseudo-differential
operators based on the asymptotic behaviour of entropy and approximation numbers, together with Carl's inequality and the Birman-Schwinger principle. The approach can be used for  the pseudo-differential operators that  factor over  a compact embedding. 

Symmetry as well as weights can be used to generate compactness of Sobolev type embeddings on $\R^d$. This was noticed   in the case of the  first order Sobolev spaces of radial functions by W.~Strauss in  the seventies of the last century, cf. \cite{Strauss}.   In the general framework of Besov and Triebel-Lizorkin spaces  a detailed study  of radial distributions has been made in \cite{SS1}, \cite{SS} and \cite{SSV}, cf. also \cite{EF1} for somewhat  different approach.   
The asymptotic behaviour of entropy numbers of the compact embedding of radial Besov spaces was described  by {Th.~K{\"u}hn}, {H.-G.~Leopold, W.~Sickel} and the second named author in \cite{KLSS}, corresponding approximation numbers  was studied in \cite{ST} and the  Gelfand and Kolmogorov numbers by the first named author in \cite{AG}.

Much less is known about  weaker symmetry assumptions, in particular about so called block-radial symmetry on $\R^d$. The compactness of the corresponding embeddings was noticed by P.L.~Lions in \cite{Lio82}. In \cite{LS} the second named author extended this result to Besov and Triebel-Lizorkin spaces. We point out that compactness of these embeddings is  a multidimensional  phenomenon, since any block should be  at least of dimension 2. Thus the simplest possible setting is the block-radial symmetry in four dimensional euclidean space with two 2-dimensional blocks.     The application of block-radial functions to nonlinear elliptic problems can be found in \cite{Ku, PK}. 

Now  our main aim  is to calculate the asymptotic  behaviour of  entropy numbers of compact embeddings of Sobolev and Besov spaces of block-radial functions and present some typical applications of this result to estimate the distribution of eigenvalues of degenerate pseudo-differential operators. Moreover   we are interested in    negative spectra of the corresponding  Schr\"odinger type operators with block-radial potentials. We estimate the number of the negative eigenvalues related to the block radial eigenfunctions.  In  particular we show that the Schr\"odinger type operators with radial potentials can have block-radial eigenvalues that are not radial. 

We recall what we mean by block-radial symmetry. 
Let $m\in\{1,\ldots ,d\}$ and let $\gamma\in\N^{m}$ be  an $m$-tuple  $\gamma=(\gamma_1,\ldots, \gamma_m)$,  $\gamma_1+\ldots + \gamma_m = |\gamma|= d$. The $m$-tuple  $\gamma$ describes the decomposition of $\R^{|\gamma|}= \R^{\gamma_{1}}\times\dots\times\R^{\gamma_{m}}$  into $m$ subspaces of dimensions $\gamma_{1},\dots,\gamma_{m}$ respectively.  Let 
\[
SO(\gamma)= SO(\gamma_1)\times \ldots \times SO(\gamma_m)\subset SO(d)
\] 
be a group of isometries on $\R^{|\gamma|}$. An element $g=(g_1,\ldots, g_m)$, $g_i\in SO(\gamma_i)$  acts  on $x=(\tilde{x}_1,\ldots, \tilde{x}_m)$, $\tilde{x}_i\in \R^{\gamma_i}$ by $x\mapsto g(x)= (g_1(\tilde{x}_1),\ldots, g_m(\tilde{x}_m))$.  If $m=1$ then $SO(\gamma)= SO(d)$ is the special orthogonal group acting on $\R^d$. If $m=d$ then the group  is trivial since then $\gamma_1= \ldots =\gamma_m = 1$ and $SO(1)= \{\mathrm{id}\}$.  We will always assume that $\gamma_i\ge 2$ for any $i=1,\ldots, m$.

Let $\bspq(\R^d)$, $s\in \R$ and $1\le p,q\le \infty$, be a Besov space and $R_\gamma\bspq(\R^d)$ be its subspaces consisted of $SO(\gamma)$-invariant distributions. It is known that if $s_1-\frac{d}{p_1}>s_2-\frac{d}{p_2}$ and $p_1<p_2$ then the embedding 
\[
R_\gamma B^{s_1}_{p_1,q_1}(\R^d) \hookrightarrow R_\gamma B^{s_2}_{p_2,q_2}(\R^d)
\] 
is compact, cf. \cite{LS}. Let us assume that  $\gamma_1\le \ldots \le\gamma_m$ and let $n= \max\{i: \gamma_i=\gamma_1 \}$. We prove  that the entropy numbers $e_k$ of this embedding have the following asymptotic behaviour
\[
e_k \, \sim \, \big( k^{-\min_i\gamma_i}(\log k)^{(n-1)(\gamma_1-1)}\big)^{\frac{1}{p_1}-\frac{1}{p_2}},
\]
cf. Theorem \ref{mainT1}.
If $m=1$ then the space $R_\gamma B^{s_1}_{p_1,q_1}(\R^d)$ consists  of radial functions and the above estimates coincide with the estimates proved in \cite{KLSS}.  Similarly to the radial case the asymptotic behaviour is independent of the smoothness parameters  $s_1$ and $s_2$.  Please note that  $\min_i \gamma_i \le \frac{d}{m}\le \max_i \gamma_i$.

In the paper we investigate also the negative spectrum of the self-adjoint operator of the type
\[
H_{s,\beta} = (\Id -\Delta)^{s/2} - \beta V\qquad\mbox{as}\quad
\beta\rightarrow\infty  .  
\] 
We show that if the potential $V\in L_r(\R^d)$ is $SO(\gamma)$-invariant and $s>\frac d r$ then the operator has asymptotically at most $\beta^{r/ \gamma_1} (\log \beta)^{(n-1)(\gamma_1-1)/\gamma_1}$ and at least $\beta^{m/s}$ negative eigenvalues with $SO(\gamma)$-invariant eigenfunctions, cf. Theorem \ref{op-main}.  If $V$ is radial and $\frac d r < s < [d/2]\frac d r$ then the operator      $H_{s,\beta,\theta}$ has eigenfunctions that are block-radial but not radial.

\subsection*{Notation}
Sobolev, Besov and Triebel-Lizorkin spaces are discussed in various places, we refer e.g. to the monographs \cite{TrI, Tr1}. We will use only  the basic definitions and facts of this theory and will not recall them here. We refer the reader to the quoted literature.    

As usual, $\N$ denotes the natural numbers, $\N_0:= \N \cup \{0\}$,  $\Z$ denotes the integers and
$\R$ the real numbers. Logarithms are always taken in base 2, $\log=\log_2$.
If $X$ and $Y$ are two Banach spaces, then the symbol  $X \hookrightarrow Y$ indicates that the embedding is continuous.
The set of all linear  and bounded operators
$T : X \to Y$, denoted by $\cl (X,Y)$, is  equipped with the standard norm.
As usual, the symbol  $c $ denotes positive constants
which depend only on the fixed parameters $s,p,q$ and probably on auxiliary functions,
unless otherwise stated; its value  may vary from line to line.
We will use the symbol $A \sim B$, where $A$ and $B$ can depend on certain parameters. The meaning of $A \sim B$ is given by: there exist  constants $c_1,c_2>0$ such that  inequalities $ c_1 A \le B \le c_2 A$ hold for all values of the parameters. 

We shall use the following conventions throughout the paper:
\begin{itemize}
	\item
	If $E$ denotes a space of distributions  (functions) on $\R^d$ then by $R_\gamma E$ we mean the subset of $SO(\gamma)$-invariant  distributions (functions) in $E$. We endow this subspace  with the same norm as the original space. If $SO(\gamma)=SO(d)$, i.e. if  the subspace consists of radial functions, then we will write $RE$.  \\
	Similarly if $G$ is a finite group of reflections in $\R^d$ then $R_G E$  denotes  the subspace of those elements of the space $E$ that are invariant with respect to $G$. 
	\item
	If an equivalence class $[f]$ (equivalence with respect to coincidence almost everywhere) contains a continuous representative
	then we call the class continuous and speak of values of $f$ at any point
	(by taking the values of the continuous representative).
	\item We will use also the following notation related to the action of the group $SO(\gamma)$ on $\R^d$, 
	\begin{align}\nonumber
	r_{j} \,=\,r_{j}(x)\, =\, \left(x_{\gamma_{1}+\dots+\gamma_{j-1}+1}^{2}+\dots+x_{\gamma_{1}+\dots+\gamma_{j-1}+\gamma_{j}}^{2}\right)^{1/2}. 
	\end{align}
\end{itemize}


\section{Traces of block-radial functions} \label{main1a}

Let $d=|\gamma|$ and $\gamma_i\ge 2$ for any $i=1,\dots ,m$.  To simplify the notation we put $\bar{\gamma}_i = 1 + \sum_{\ell=0}^{i-1} \gamma_\ell$ if $i=1,2,\ldots, m+1$ with ${\gamma}_0= 0$. We define a hyperplane 
\begin{equation}\nonumber
H_\gamma = \mathrm{span}\{e_{\bar{\gamma}_1},e_{\bar{\gamma}_2}, , \ldots , e_{\bar{\gamma}_m}\}\ ,
\end{equation}
where $e_j$, $j=1,\ldots, d$ is a standard orthonormal basis in $\R^d$. The hyperplane $H_\gamma$ can be  identified   with $\R^m$ in the standard way so we  write $(r_1,\ldots , r_m)\in H_\gamma$ if $r_1 e_{\bar{\gamma}_1} + r_2 e_{\bar{\gamma}_2} + \ldots + r_m e_{\bar{\gamma}_m}\in H_\gamma$.  We  need also a finite group of reflections  $G(\gamma)$ acting on $H_\gamma$. The group consists of  transformations  $g_{i_1,\ldots,i_m}\in G(\gamma)$ given by  \[
g_{i_1,\ldots,i_m}(r_1,\ldots , r_m) = \Big( (-1)^{i_1}r_1,\dots , (-1)^{i_m}r_m\Big), \qquad  (i_1,\ldots,i_m)\in \{1,2\}^{m}\, . 
\]

Let $f:\Rd \to \C$ be a locally integrable $SO(\gamma)$-invariant function.
By using the Lebesgue point argument its restriction
\[
f_0 (r_1,\ldots, r_m):= f(\widetilde{r}_1,\ldots, \widetilde{r}_m) \ , \qquad  \widetilde{r}_j=(r_j,0, \ldots \, , 0)\in \R^{\gamma_j}\ , \quad  j=1,\ldots , m\, 
\]
is well-defined a.e. on $H_\gamma$.
However, this restriction need not be locally integrable.
A simple example is given by the function
\[
f(x):= \psi (x)\, |x|^{-m}\, , \quad x \in \Rd , \quad \psi\in C_o(\R^d), \quad \psi(0)=1\ .
\]
On the other hand, we can start with a measurable  $g:\, H_\gamma \to \C$, that is invariant with respect to the action of the group $G(\gamma)$. 
If  $g$ is locally integrable on all subsets $\{(r_1,\ldots ,r_m): \; r_j>0,\; j=1,\ldots ,m, \; \text{and}\;  a< |(r_1,\ldots ,r_m)| <b\;\}$, $0 < a < b < \infty$,
then (again using the Lebesgue point argument) the function
\[
f(x):= g(r(x)) \, , \qquad x \in \Rd\,
\]
is well-defined a.e. on $\Rd$ and is  $SO(\gamma)$-invariant.
In what follows we shall study properties of the associated  operators
\[
\tr \, : f \mapsto f_0 \qquad \mbox{and}\qquad \ext \, : \, g \mapsto f\, .
\]
Both operators are defined pointwise.


\subsection{Traces of block-radial $L_p$-spaces.}
\label{mainL_p}


Before  we turn to the description of the trace classes of block-radial Besov and Sobolev 
spaces with $1\leq p\leq \infty$ we start with  almost trivial results for $L_p$-functions.  We need a further notation.
By $L_p (\R^m,w)$ we denote the weighted Lebesgue space equipped with the norm
\[
\|\, f \, |L_p (\R^m,w)\| := \Big(\int_{\R^m} |f(x)|^p\, w(x)\, dx\Big)^{1/p}
\]
with the usual modification if $p=\infty$.
We will use  a weight
\begin{equation} 
w_\gamma(r_1,\dots , r_m)\  = \ \prod_{i=1}^{m}|r_i|^{\gamma_i-1}.
\label{wg}
\end{equation}
Direct calculations show that $w_\gamma$ is a Muckenhoupt weight. More precisely $w_\gamma\in {\mathcal A}_\rho$ for any $\rho>\max \gamma_i$. We recall the definition of  the  ${\mathcal A}_\rho$ classes in Appendix B.

\begin{Lem}\label{simple}
	We assume that $d\ge 2$.\\
	{\rm (i)}
	Let $0< p< \infty$. Then
	$\tr \, :\,  R_{\gamma}L_p (\Rd) \to R_{G}L_p (H_\gamma, w_\gamma)$ is a linear isomorphism with
	inverse  $\ext$.
	\\
	{\rm (ii)}
	Let $p = \infty$. Then $\tr \, :\,  R_{\gamma}L_\infty (\Rd) \to R_{G}L_\infty (H_\gamma)$ is a linear isomorphism with  inverse $\ext $.
\end{Lem}

\bpr
\noindent
Introducing the radial coordinates on each block we get 
\begin{align}\label{Lptrace}
\int_{\Rd} &|f(x)|^p \, dx  =   \\ & \frac{2\pi^{d/2}}{\Gamma(\gamma_1/2)\ldots \Gamma(\gamma_m/2)} \, \int_0^\infty\ldots \int_0^\infty |f_0(r_1,\ldots,r_m)|^p \, r_1^{\gamma_1-1}\ldots r_m^{\gamma_m-1}\, dr_1\ldots dr_m\, .\nonumber
\end{align}
On the other hand the formula 
\begin{align*}
\int_0^\infty\ldots \int_0^\infty & |f_0(r_1,\ldots,r_m)|^p \, r_1^{\gamma_1-1}\ldots r_m^{\gamma_m-1}\, dr_1\ldots dr_m\, =\\  
& = \lim_{\varepsilon \downarrow 0} \int_\varepsilon^\infty\ldots \int_\varepsilon^\infty |f_0(r_1,\ldots,r_m)|^p \, r_1^{\gamma_1-1}\ldots r_m^{\gamma_m-1}\, dr_1\ldots dr_m\, ,
\end{align*}
implies that   test functions  supported in the interior of $[0, \infty)^m$ are dense in $ L_p([0, \infty)^m, w_\gamma)$. In consequence 
the  formula \eqref{Lptrace} can be read  from the other side i.e.
\begin{align*}
\int_{\Rd} & | \ext g (x)|^p \, dx  = \\
& =  \frac{2\pi^{d/2}}{\Gamma(\gamma_1/2)\ldots \Gamma(\gamma_m/2)} \, \int_0^\infty\ldots \int_0^\infty |g(r_1,\ldots,r_m)|^p \, r_1^{\gamma_1-1}\ldots r_m^{\gamma_m-1}\, dr_1\ldots dr_m\, ,
\end{align*}
for all $g \in  L_p ([0, \infty)^m, w_\gamma )$. This proves (i).
Part (ii) is obvious.
\epr

\noindent
Lemma \ref{simple} means that whenever  the Besov $\bspq (\Rd)$ or Sobolev  space $H^s_{p}(\Rd)$
is contained in $ L_1  (\Rd) +L_\infty (\Rd)$,
then $\tr $ is well-defined on its block-radial subspace.
It is well known that 
\[
\bspq (\Rd)\, , \, H^s_p(\Rd) \hookrightarrow  L_1  (\Rd) +L_\infty (\Rd)
\]
if $s>d \, \max(0,\frac 1p -1)$, see e.g. \cite{SiTr}.

This is in some contrast to the general theory of traces on
these spaces. Generally to guarantee that $\bspq (\Rd)$ or $H^s_p(\Rd)$ has a trace  on $\R^m$ one has to assume that 
\[
s > \frac{d-m}{p} + m \max \Big(0, \frac 1p -1\Big)\, ,
\]
cf. e.g. 
\cite{FJ1},   \cite[Rem.~2.7.2/4]{Tr1}. 
This condition is stronger.

\subsection{Traces of block-radial Sobolev and Besov spaces} 

Let $n \in \N_0$. Then 
$W^n_p (\Rd)$ denotes the collection of all functions  $f: \Rd \to \C$
such that all weak derivatives $D^\alpha f$ of order $|\alpha| \le n$ exist
and belong to $L_p(\R^d)$. 
The norm in $W^n_p (\Rd)$ is defined by 
\begin{equation} \nonumber \label{sobnorm}
\| \, f \, |W^n_p(\Rd)\| := \sum_{|\alpha|\le n} \| \, D^\alpha f\, |L_p (\Rd)\|\, .
\end{equation}

\begin{Prop}
	\label{anfang}
	Let $d\ge 2$.
	For $1<p<\infty$ and $n >  \max(0, \frac{\max \gamma_i}{p}-1)$, $n  \in \N$  the  mapping
	$\tr$ is a linear  isomorphism of $R_{\gamma}W^n_p (\Rd)$ onto $R_{G}W^n_p (H_\gamma, w_\gamma)$
	with the  inverse  $\ext$.
\end{Prop}


\bpr
{\em Step 1.} (The trace operator.) 
The operator $\tr$ is well-defined on  $R_{\gamma}W^n_p (\Rd)$ and $\tr\big(R_{\gamma}W^n_p (\Rd)\big)\subset R_{G}L_p (H_\gamma, w_\gamma)$ since $R_{\gamma}W^n_p (\Rd)\subset  R_{\gamma}L_p (\Rd)$, cf. Lemma \ref{simple}. 
One can easily see that  if $f\in R_{\gamma}C^n_0 (\Rd)$ then $f_0=\tr\, f \in R_{G}C^n_0 (H_\gamma)$ since 
$$
\frac{\partial^{|\alpha|}f_0}{\partial r_m^{\alpha_m}\ldots  \partial r_1^{\alpha_1}} (r(x)) \, =\,  
\frac{\partial^{|\alpha|}f}{\partial x_{\bar{\gamma}_m}^{\alpha_m}\ldots  \partial x_{\bar{\gamma}_1}^{\alpha_m}} (x) \ , \qquad |\alpha|\le n\ . 
$$
In consequence introducing the block spherical coordinates \eqref{Lptrace} we get 
\begin{equation}\label{13.11_0}\nonumber
\|f_0 |W^n_p(\R^m, w_\gamma)  \| 
\, \le \, C\, \|f | W^n_p(\R^d) \|\,. 
\end{equation}
It was proved in \cite{LS} that the space $R_{\gamma}W^n_p(\R^d)$ is a complemented subspace of $W^n_p(\R^d)$ and that the corresponding projection maps $C^\infty_0(\R^d)$  onto  $R_{\gamma}C^\infty_0(\R^d)$. Thus the space $R_{\gamma}C^\infty_0(\R^d)$ is dense in $R_{\gamma}W^n_p(\R^d)$.  This proves that 
$$
\tr: \ R_{\gamma}W^n_p (\Rd) \ \rightarrow R_{G}W^n_p (H_\gamma, w_\gamma)\, . 
$$ 

{\em Step 2.} (The extension operator)
Since $1<p<\infty$ one can define an equivalent norm in the Sobolev spaces by
\begin{equation}\label{sobnorm1}
\| \, f \, |W^n_p(\Rd)\|_{(1)} := \| \, f\, |L_p (\Rd)\| + \sum_{i=1}^d \left\| \, \frac{\partial^n f}{\partial x^n_i}\, |L_p (\Rd)\right\|\, .
\end{equation}

Let $f_0\in  R_{G}C^n_0 (H_\gamma)$ and $ f = \ext f_0$. First we  fix $(x_{\gamma_1+1}, \ldots, x_d)\in \R^{d-\gamma_1}$. The function
\begin{equation}\label{1111-1}\nonumber
\widetilde{f} (x_1,\ldots , x_{\gamma_1}): = f((x_1,\ldots , x_{d})) = f_0(r_{1}(x), \ldots , r_m(x)) 
\end{equation} 
is radial on $\R^{\gamma_1}$.  Moreover $\widetilde{f}$ is the radial extension of an even function $\widetilde{f_0}(r)= f_0(r, r_2(x),\ldots , r_m(x))$, $r\in \R$.  
By the trace result for radial functions, Theorem 3, Theorem 8 and Theorem 9 in \cite{SSV}, we  deduce  that there is a constant $C>0$ independent of $f_0$ such that the inequality 
\begin{align}\label{13.11_1}
& \int_{\R^{\gamma_1}} \left|\partial^\alpha f (x_1,\ldots ,x_d)\right|^p d x_1\ldots d x_{\gamma_1} \, \le \, C\, 
\|\widetilde{f_0}| W^n_p(\R, |r|^{\gamma_1-1} )\|^p = \\  \qquad \qquad & =  C \sum_{i=0}^n \int_{\R} \left|\frac{d^i}{d r^i} f_0(r, r_2(x),\ldots ,r_m(x))\right|^p |r|^{\gamma_1-1}d r   \nonumber
\end{align} 
holds for any $\alpha\in \{1,\ldots , \gamma_1\}^{\gamma_1}$, $|\alpha|\le n$. The function 
\[
\R^d\ni x \mapsto f_0(r, r_2(x),\ldots ,r_m(x))
\]
is invariant with respect of any isometry belonging to  $\{\id \}\times SO(\gamma_2)\times \ldots \times SO(\gamma_m)$, therefore integrating the inequality  \eqref{13.11_1} with respect to the variables $(x_{\gamma_1+1}, \ldots, x_d)$  we obtain
\begin{align}\label{13.11_1a}\nonumber
& \int_{\R^{d}} \left|\partial^\alpha f (x_1,\ldots ,x_d)\right|^p d x_1\ldots d x_{d} \, \le  \\  
& \qquad  \le   C \sum_{i=0}^n \int_{\R} \left(\int_{\R^{d-\gamma_1}} \left|\frac{d^i}{d r^i} f_0(r, r_2(x),\ldots ,r_m(x))\right|^p d x_{\gamma_1+1}\ldots d x_{d} \right) |r|^{\gamma_1-1}d r  =  \nonumber \\   
& \qquad    = C  \sum_{i=0}^n \int_{\R^m}  \left|\frac{d^i}{d r_1^i} f_0(r_1, r_2,\ldots ,r_m)\right|^p |r_1|^{\gamma_1-1}\ldots |r_m|^{\gamma_m-1} d r_1\ldots d r_m  \,\le  \nonumber \\  
& \qquad \qquad  \le \, C\, \|f_0| W^n_p(\R^m, w_\gamma)  \|^p . \nonumber
\end{align} 

Similar argument works for any $\alpha\in \{\bar{\gamma}_i,\ldots , \bar{\gamma}_{i+1}-1\}^{\gamma_i}$, $i=2,\ldots , m-1$ , $|\alpha|\le n$. So for any $j=1,\ldots , d$ we have
\begin{equation}\nonumber
\int_{\R^{d}} \left|\frac{\partial^n}{\partial x_j^n} f (x_1,\ldots ,x_d)\right|^p d x_1\ldots d x_{d} \, \le  C\, \|f_0 | W^n_p(\R^m, w_\gamma)  \|^p.
\end{equation} 

Now using the norm \eqref{sobnorm1} we conclude
\begin{equation}\label{13.11_3}\nonumber
\|\ext f_0 | W^n_p(\R^d)  \| = \|f | W^n_p(\R^d)  \| \, \le \, C\, \|f_0 | W^n_p(\R^m, w_\gamma)\|\,.  
\end{equation}

Smooth compactly supported functions are dense in $W^n_p(\R^m, w_\gamma)$, cf. \cite{Bui}. The space  $R_{G}W^n_p(H_\gamma, w_\gamma)$ is a complemented subspace of $W^n_p(H_\gamma, w_\gamma)$ since 
\[
P: W^n_p(\R^m, w_\gamma)\ni h \mapsto \frac{1}{|G(\gamma)|}\sum_{g\in G(\gamma)} h\circ g  
\]
is a continuous projection onto $R_{G}W^n_p(H_\gamma, w_\gamma)$. Moreover $P$ maps $C^\infty_0(\R^m)$ onto $R_{G}C^\infty_0(H_\gamma)$. Thus $R_{G}C^\infty_0(H_\gamma)$ is a dense subspace of $R_{G}W^n_p(H_\gamma, w_\gamma)$. Thus we can  extend the operator $\ext$ to a continuous operator defined on $R_{G}W^n_p(H_\gamma, w_\gamma)$. This proves the proposition. \epr

\begin{Cor}
	\label{anfang1}
	Let $d\ge 2$. Let $1<p<\infty$,  $s>  \big[\max(0, \frac{\max \gamma_i}{p}-1)\Big] +1$, and $0<q\leq \infty$.  
	
	\noindent
	(a) The  mapping $\tr$ is a linear  isomorphism of $R_{\gamma} H^s_p (\Rd)$ onto $R_{G}H^s_p (H_\gamma,  w_\gamma)$
	with  the inverse  $\ext$.
	
	\noindent
	(b) The  mapping $\tr$ is a linear  isomorphism of $R_{\gamma} B^s_{p,q}(\Rd)$ onto $R_{G} B^s_{p,q} (H_\gamma,  w_\gamma)$
	with  the inverse  $\ext$.
\end{Cor}

\bpr
The point (a) can be proved by complex interpolation  and the retraction-coretraction method. 
We known that 
\begin{equation}\nonumber
\left[ W^{n_1}_p(\R^m, w), W^{n_2}_p(\R^m, w)\right]_\theta = H^s_p(\R^m, w)\qquad\text{if}\quad (1-\theta) n_1 + \theta n_2  = s ,  
\end{equation}  
with $w$ being the Muckenhoupt weight, cf. \cite{SSV2}. Moreover it was proved in \cite{LS} that the space $R_{\gamma}W^n_p(\R^d)$ is a complemented subspace of $W^n_p(\R^d)$. Analogously $R_{G}W^n_p(H_\gamma,w_\gamma)$ is a complemented subspace of $W^n_p(H_\gamma,w_\gamma)$. So using the retraction-coretraction method we can prove that  
\begin{align*}
\left[ R_{\gamma}W^{n_1}_p(\R^d), R_{\gamma}W^{n_2}_p(\R^d)\right]_\theta & = R_{\gamma}H^s_p(\R^d)\qquad\text{if}\quad (1-\theta) n_1 + \theta n_2  = s \, , \\     
\left[ R_{G}W^{n_1}_p(H_\gamma, w_\gamma), R_{G}W^{n_2}_p(H_\gamma, w_\gamma)\right]_\theta  & = R_{G}H^s_p(H_\gamma, w_\gamma)\qquad\text{if}\quad (1-\theta) n_1 + \theta n_2  = s \, .     
\end{align*}  

The point (b) can be proved similarly, but know we use the real interpolation. More precisely one should use the formula 
\[
\left( W^{n_1}_p(\R^m, w), W^{n_2}_p(\R^m, w)\right)_{\theta, q} = B^s_{p,q}(\R^m,  w)\qquad\text{if}\quad \theta k = s , \; k>s \, , \; 0<q\le \infty ,     
\] 
proved by Bui in \cite{Bui}. \epr

Now we consider the case $p=\infty$.

\begin{Prop}
	\label{anfang2}
	Let $d\ge 2$, $0<s<1$, $p = \infty$ and $0<q \le \infty$.  The  mapping
	$\tr$ is a linear  isomorphism of $R_{\gamma} B^s_{\infty,q}(\Rd)$ onto $R_{G} B^s_{\infty,q} (H_\gamma,  w_\gamma)$ with  the inverse  $\ext$.
\end{Prop}

\bpr   
It is obvious that $\tr$ is a linear  isomorphism of $R_{\gamma}C (\Rd)$ onto $R_{G}C (H_\gamma)$ with  the inverse  $\ext$.
We prove that it is also a linear isomorphism of  $R_{\gamma}C^1 (\Rd)$ onto $R_{G}C^1 (H_\gamma)$ with  the same inverse. Here  
$C^1 (\Rd)$ denotes the collection of all functions  $f: \Rd \to \C$ such that all their derivatives of  the first order  exist, are uniformly continuous and bounded. The spaces $C^1 (\Rd)$ is equipped with the norm 
\[
\| \, f \, |C^1(\Rd)\| := \sum_{|\alpha|\le 1} \| \, D^\alpha f\,|L_\infty (\Rd)\|\, .
\]

It should be clear that if $f$ is $SO(\gamma)$-invariant on $\Rd$ then $\tr f$ is a continuous $G$-invariant function on  $H_\gamma$. Vice versa,  if $f_0\in C^1 (H_\gamma)$ is $G$-invariant on  $H_\gamma$ then $\ext f_0$ is a continuous $SO(\gamma)$-invariant on $\Rd$.

Let $f \in R_{\gamma} C^1 (\Rd)$.
For $r=(r_1,\dots ,r_m)\in H_\gamma$  we put $\tilde{r}=(x_1,\dots , x_d)$ with $x_\ell= r_j$ if $\ell=\bar{\gamma}_j$ and $x_\ell=0$ otherwise. 
We obviously have
\[
\frac{\partial f_0}{\partial r_j}(r_1,\dots ,r_m) = \frac{\partial f}{\partial x_{\bar{\gamma}_j}} (\tilde{r})  \ ,
\]
which proves the estimate
\begin{equation} \label{eq:trf} \nonumber
\|\, \tr f \, |C^1 (H_\gamma)\| \le \| \, f \, |C^1 (\R^d)\|
\end{equation}
and at the same time the continuity of the function  $\tr f= f_0$ and its derivative.\\
Now, we assume that $f_0 \in R_{G}C^1 (H_\gamma)$. Let $f:= \ext f_0$.
If $r_j(x) \neq 0$ then we have
\be\label{eq-41b}
\frac{\partial f}{\partial x_\ell} (x) =
\frac{\partial f_0}{\partial r_j}(r(x))\,  \frac{x_\ell}{r_j(x)} \, , \qquad \bar{\gamma}_j\le \ell < \bar{\gamma}_{j+1}\, .
\ee

\noindent
Let $r_j(x) = 0$. The function $f_0$ is $G$-invariant,  so the function $r_i\mapsto f_0(r_1,\dots ,r_m)$ is even for any $i=1,\ldots , m$. In consequence  derivatives of the  continuously differentiable function $f_0$ satisfies 
\begin{equation}\label{Cderiv1}
\frac{\partial f_0}{\partial r_j}(r_1,\ldots ,r_{j-1},0,r_{j+1}, \ldots, r_m) = 0\, .  
\end{equation}
So, if $\bar{\gamma}_j\le \ell < \bar{\gamma}_{j+1}$ then \eqref{Cderiv1} implies
\begin{align}\label{Cderiv2}
&\frac{\partial f}{\partial x_\ell}(x)  =  \lim_{h\to 0} \frac{f(x_1,\ldots ,x_\ell+h, \ldots, x_d) - f(x_1,\ldots , x_d)}{h} =\\
& \qquad = \lim_{h\to 0} \frac{f_0(r_1,\ldots ,r_{j-1},|h|,r_{j+1}, \ldots, r_m) - f_0(r_1,\ldots ,r_{j-1},0,r_{j+1}, \ldots, r_m)}{h} = 0. \nonumber 
\end{align}

Now \eqref{eq-41b} and \eqref{Cderiv2} give us 
\[
\sup_{x }\, \Big|\frac{\partial f}{\partial x_\ell} (x)\Big|\le
\sum_{j=1}^{m}\sup_{r}\, \Big|\frac{\partial f_0}{\partial r_j}(r)\Big| \, .
\]

It remains to deal with the continuity of the derivatives $\frac{\partial f}{\partial x_\ell}$. Let   $\bar{\gamma}_j\le \ell < \bar{\gamma}_{j+1}$. If $r_j(x)\not= 0$ then the continuity follows immediately from \eqref{eq-41b} and the continuity of derivatives of function $f_0$. 

\noindent 
Let $r_j(x)= 0$ and let  $x^{(k)}\rightarrow x$ in $\R^d$ as $k\rightarrow \infty$, $r_j(x^{(k)})\not= 0$.  Then $r_j(x^{(k)})\rightarrow 0$ so \eqref{eq-41b}, \eqref{Cderiv1} and the continuity of the partial derivatives of $f_0$ imply the continuity of the partial derivative of $f$ at $x$.   
This proves the claim.

It remains to extend the statement to Besov spaces $B^s_{\infty,q}$.  Once more this can be done via interpolation since
\[ 
(C(\R^d), C^1(\R^d))_{s, q} = B^s_{\infty, q} (\R^d), \qquad 0<s<1 \qquad\text{and}\quad 0<q\le \infty \, , 
\]
cf. Theorem 2.7.2 and Theorem 1.10.2 in  \cite{TrI} . 
\epr

\begin{Cor}\label{anfang3}
	Let $d\ge 2$, $s\in \R$, $1 < p \le \infty$ and $0<q \le \infty$.  Then the space 
	$R_{\gamma} B^s_{p,q}(\Rd)$ is isomorphic to  $R_{G} B^s_{p,q} (H_\gamma,  w_\gamma)$.  Analogously the space  
	$R_{\gamma} H^s_p (\Rd)$ is isomorphic to $R_{G}H^s_p (H_\gamma,  w_\gamma)$. 
\end{Cor}

\bpr
If $1<p<\infty$ and $s$ is sufficiently large  then the corollary  follows immediately form Corollary \ref{anfang1}.  For other  value of $s$ the statement follows by the lift property. 
The statement can be extented to smaller values of $s$  by  the lifting property for (weighted) Besov and Sobolev  spaces, cf. \cite{Bui}. 

The same argument works for $p=\infty$. Now we should use Proposition \ref{anfang3} instead of  Corollary \ref{anfang1}.   
\epr


\section{Entropy  numbers  of embeddings of spaces of block-radial functions}

In this section we estimate   entropy  numbers  of compact Sobolev  embeddings of block-radial functions. It was proved in \cite{LS} that the embedding    
\begin{align}
\label{ent1}
& id~:~R_{\gamma} B^{s_1}_{p_1,q_1}(\Rd)\rightarrow R_{\gamma} B^{s_2}_{p_2,q_2}(\Rd) 
\intertext{is compact if and only if} 
& p_1< p_2,   \quad\delta = s_1- \frac{d}{p_1} - s_2+ \frac{d}{p_2}>0 \quad \text{and}\quad \min_i\gamma_i\ge 2\, . \label{ent1a}
\end{align}

For convenience  of the reader we recall the basic definitions. Let $X$ and $Y$ be Banach spaces and $T \in {\cal L}(X,Y)$. The $k$-th entropy number of $T$, $k\in\N$, is defined in the following way 
$$ 
e_k(T):=\inf\{\epsilon>0 : T(B_X) \text{\ can be covered by\ } 2^{k-1} \text{\ balls of radius\ } \epsilon \text{\ in\ } Y\},
$$
where $B_X$ denotes the closed unit ball in $X$. The sequence of entropy numbers can be viewed as  quantification of the notion of compactness since the operator is compact if and only if  $ e_k(T) \rightarrow 0,$ as $ k \rightarrow \infty$.

The entropy  numbers have properties of   multiplicativity and additivity, i.e. 
\begin{align}\nonumber
e_{n+k-1}(T\circ S) & \le e_n(T) \cdot e_k(S) \qquad\text{and}\qquad 
e_{n+k-1}(T+S) & \le e_n(T) + e_k(S) .
\end{align}
For further  their properties  we refer to \cite{Pi}, \cite{Koe}, \cite{CS} or \cite{ET}.

To estimate entropy numbers we use the technique of quasi-normed operator ideals. In the context of entropy numbers the approach goes back to  B. Carl's paper \cite{Carl2}.  For Sobolev embeddings it was used in \cite{Kuh3} and \cite{KLSS} for the first time.   
Let $ \omega = (\omega_n) $ be an
increasing sequence of positive real numbers satisfying the
regularity condition $\omega_{2k}\sim\omega_k$. The properties of entropy numbers  imply that for operators
$T\in {\cal L}(X,Y)$ between Banach spaces the formula 
$$
L^{(e)}_\omega(T)\,:=\,  \sup_k \,\omega_k e_k(T)\,
$$
defines a quasi-norm in the vector space
$$
{\cal L}^{(e)}_\omega(X,Y)=\{\,T\in{\cal L}(X,Y)\,: L_\omega
^{(e)}(T)<\infty\,\}\,.
$$ In fact, $\big({\cal
L}^{(e)}_\omega,L^{(e)}_\omega\big)$ is a quasi-normed operator
ideal in the sense of Pietsch \cite[Definition 6.1.1]{Pi2}.

By Corollary \ref{anfang3} we can reduce investigation  of an asymptotic behaviour of entropy numbers of the embeddings \eqref{ent1} to  estimation of embeddings of the  corresponding weighted spaces with the Muckenhoupt weight $w_\gamma$, cf. \eqref{wg},   i.e.  
\begin{align*}
e_k(id: R_{\gamma} B^{s_1}_{p_1,q_1}(\Rd)&\rightarrow R_{\gamma} B^{s_2}_{p_2,q_2}(\Rd))\sim \\ & \sim e_k(id:R_G B^{s_1}_{p_1,q_1}(\R^m, w_\gamma)\rightarrow R_G B^{s_2}_{p_2,q_2}(\R^m, w_\gamma)).
\nonumber
\end{align*}
Furthermore using the wavelet characterization of Besov spaces with ${\mathcal A}_\infty$ weights we can use the technique of discretization i.e., we can reduce the problem to the corresponding problem for suitable sequence spaces,  cf. 
\cite[Theorem 1.13]{HS}. However, the resulting sequence spaces are still complicated, therefore a further reduction is necessary. We will use the following result concerning the entropy numbers of general diagonal operators proved by Th.~K\"uhn,  cf. \cite{Kuh, Kuh2}.      

\begin{T}[cf. \cite{Kuh2}]\label{Kuhn}
Let $0< p_1,p_2\leq \infty$, and let $\sigma=(\sigma_k)$  be a non-increasing sequence satisfying the doubling condition $\sigma_k\sim \sigma_{2k}$ and, in addition,
\begin{equation}
\sup_{n\geq k}\frac{\sigma_n}{\sigma_k}\cdot\big(\frac n k \big)^\alpha<\infty\quad\text{for some}\ \alpha>\max(1/p_2-1/p_1,0).
\label{adw}
\end{equation}
Then 
	$$
	e_k(D_\sigma: \ell_{p_1}\rightarrow\ell_{p_2})\sim k^{\frac{1}{p_2}-\frac{1}{p_1}}\sigma_k.
	$$
\end{T}

Now we recall the definition of the sequence spaces.    Let   $Q_{\nu,n}$  denote a dyadic cube in $\R^m$, centred at $2^{\nu}n$, $n\in\Z^m$, $\nu\in\N_0$,  and with the side length $2^{-\nu}$.  For $0<p<\infty$, $\nu\in\N_0$ and $n\in\Z^m$ we denote by $\chi_{\nu,n}^{(p)}$ the $p$-normalized characteristic function of the cube $Q_{\nu,n}$,
\[ 
\chi_{\nu,n}^{(p)}(x)=2^{\frac{m \nu}{p}}\chi_{\nu,n}(x)=\begin{cases}
2^{\frac{m\nu }{p}}, &\text{for} \quad x\in Q_{\nu,n},\\  
0, &\text{for}\quad x\notin Q_{\nu,n}.  
\end{cases}
\]
For $0<p<\infty$, $0<q\leq\infty$, $\sigma\in\R$ and $w\in {\mathcal A}_\infty$ we introduce  sequence spaces $b_{p,q}^\sigma(w)$  
by
\begin{align*}
b_{p,q}^\sigma(w):=\biggl{\{}& \lambda=\{\lambda_{\nu,n}\}_{v,n}: \lambda_{\nu,n}\in\C,\\
&\|\lambda|b_{p,q}^\sigma(w)\|=\bigg(\sum_{\nu=0}^\infty 2^{\nu\sigma q}\bigg\|\sum_{n\in\Z^m} \lambda_{\nu,n}\chi_{\nu,n}^{(p)}\bigg|L_p(\R^m,w)\bigg\|^q\bigg)^{1/q}<\infty\biggr{\}}.
\end{align*}
If $\sigma=0$ we write $b_{p,q}(w)$ instead of $b_{p,q}^\sigma(w)$; moreover, if $w\equiv 1$ we write $b_{p,q}^\sigma$ instead of $b_{p,q}^\sigma(w)$.
Using the same arguments as in \cite{HS} we can prove that
\begin{align}\label{Bb}\nonumber
e_k(id: B^{s_1}_{p_1,q_1}(\R^m, & w_\gamma)\rightarrow  B^{s_2}_{p_2,q_2}( \R^m, w_\gamma)) \sim  e_k(id: b^{\sigma_1}_{p_1,q_1}(\ w_\gamma)\rightarrow  b^{\sigma_2}_{p_2,q_2}( w_\gamma) ), 
\end{align}
where $\sigma_i= s_i +\frac{m}{2}- \frac{m}{p_i}$, $i=1,2$.  

For later use we introduce an abbreviation
\begin{equation}\nonumber
w_\gamma(\Omega) = \int_\Omega w_\gamma(x) dx,
\end{equation}
where $\Omega\subset\R^m$ is some bounded, measurable set. For any dyadic cube $Q_{\nu,n}$ we have
\begin{align}\label{wQ}
w_\gamma(Q_{\nu,n}) &\sim  \prod_{i=1}^m \int_{2^{-\nu} n_i}^{2^{-\nu}(n_i+1)}  |r_i|^{\gamma_i-1}d r_i\sim  \prod_{i=1}^m\int_{0}^{2^{-\nu}} (r_i+2^{-\nu} n_i)^{\gamma_i-1}d r_i  \\ 
&\sim 2^{-\nu (\gamma_1+\ldots +\gamma_m)} \prod_{i=1}^m \max\{1,|n_i|\}^{\gamma_i-1} \, \sim \, 
2^{-\nu d} w_\gamma(Q_{0,n}) .\nonumber
\end{align}
Moreover, one can easily verify that the expression
\begin{align}
\label{norm}\nonumber
\bigg(\sum_{\nu=0}^\infty 2^{\nu\sigma q}\bigg(\sum_{n\in\Z^m} |\lambda_{\nu,n}|^p2^{m\nu}w_\gamma(Q_{\nu,n})\bigg)^{q/p}\bigg)^{1/q}
\end{align}
is an equivalent (quasi)-norm in $b_{p,q}^\sigma(w)$. We will used this norm in the sequel.  Please note that the conditions \eqref{ent1a} are equivalent to the necessary and sufficient conditions for  compactness of the embedding 
$id: B^{s_1}_{p_1,q_1}(\R^m, w_\gamma)\rightarrow  B^{s_2}_{p_2,q_2}(\R^m, w_\gamma)$ given in \cite[Proposition 3.1]{HS}. 

\begin{Lem}\label{volumn}
Let $\gamma= (\gamma_1 , \gamma_2, \ldots , \gamma_m)\in \N^m$ be a multi-index such that $\gamma_i\ge 2$ for any $i=1,\ldots ,m$.  
We assume that   
\begin{equation}\label{gamma-ineq}
\gamma_1 \le  \gamma_2\le  \ldots \le  \gamma_m
\end{equation}
and put 
\begin{equation}\label{gamma-ineq-n}
n = \max\{ i\,:\, \gamma_i=\gamma_1\}.
\end{equation} 
  
Let  $\tau$ denote a bijection of $\Z^m$ onto $\N_0$ such that $\tau(k)< \tau(\ell)$ if $w_\gamma (Q_{0,k}) <  w_\gamma (Q_{0,\ell})$. 
Then  there are positive constants $c_1$ and $c_2$ such that for sufficiently large $L\in \N$ the inequalities 
\begin{equation}\label{red1}
c_1 2^{L(d-m)} \le w_\gamma (Q_{0,k}) \,\le \,c_2 2^{L(d-m)}  
\end{equation} 
 if and only if 
\begin{equation}\label{red1a}
c_1  2^{\frac{L(d-m)}{\gamma_1-1}}L^{n-1}\le \tau(k)  \,\le \,c_2 2^{\frac{L(d-m)}{\gamma_1-1}}L^{n-1}.
\end{equation}
\end{Lem}

\bpr
\emph{Step 1.} To simplify our notation we put  $\alpha_i=\frac{\gamma_i-1}{d-m}$. By assumptions $0<\alpha_1\le \alpha_2\le \ldots \le \alpha_m<1$. We consider two sets 
\begin{align*}
W_L = \, & \{ (x_1,\ldots ,x_m)\in \R^m: \quad  2^L \le \prod_{i=1}^m \max(1,|x_i|)^{\alpha_i}< 2^{L+1} \} , 
\intertext{and}
\widetilde{W}_L = \, & \{ (x_1,\ldots ,x_m)\in \R^m: |x_i|\ge 1\quad\text{and}\quad  2^L \le \prod_{i=1}^m |x_i|^{{\alpha_i}}< 2^{L+1} \} .
\end{align*}
If $(k_1,\ldots, k_m)\in \N_0^m$ belongs to   $W_L$ then  $(\max(1,k_1),\ldots, \max(1,k_m))\in \widetilde{W}_L$, in consequence  
\begin{equation}\label{red2}\nonumber
\#\{k\in \N_0^m: k\in \widetilde{W}_L \} \le \#\{k\in \N_0^m: k\in W_L \} \le 2^m\#\{k\in \N_0^m: k\in \widetilde{W}_L \}.
\end{equation} 
We prove that
\begin{equation}\label{red3}
\#\{k\in \N_0^m: k\in \widetilde{W}_L \} \sim 2^{L/{\alpha_1}}L^n .
\end{equation}

First let us note that  
\begin{equation}\label{red3a}
\#\{k\in \N_0^m: k\in \widetilde{W}_L \}\sim \vol_m(\widetilde{W}_L ) = \vol_m(V^m_{2^{L+1}}) - \vol_m(V^m_{2^L}),
\end{equation}
where $V^m_{2^L}=\{ (x_1,\ldots ,x_m)\in \R^m: |x_i|\ge 1\quad\text{and}\quad  \prod_{i=1}^m |x_i|^{{\alpha_i}}\le 2^{L} \}$.

It is sufficient to restrict our attention to the first octant $\{(x_1,\ldots, x_m): x_i\ge 0, \; i=1,\ldots m \} $.    A bit more generally we consider the sets 
$$V^n_R=\{ (x_1,\ldots ,x_n)\in \R^n: x_i\ge 1\quad\text{and}\quad  \prod_{i=1}^n x_i^{{\alpha_i}}\le R \},\quad  R>1,\quad n=1,\ldots ,m ,$$ 
and show that 
\begin{equation}\label{VmR}
\vol_m(V^m_R)\sim R^{1/\alpha_1}(\log R)^{n-1}, 
\end{equation}
for sufficiently large $R$. 

\emph{Step 2.} First we consider the special case $n=m$ i.e. $\alpha_1=\ldots =\alpha_m$. The volume estimate in this case has been already calculated by Th. K\"uhn,W.Sickel and T.Ullrich in  \cite{KSU}, cf. Lemma 3.2 ibidem. Let 
\[
\mathcal{H}_\ell(r):=\{x\in [1,\infty)^\ell:  \prod_{i=1}^\ell x_i\le r \} ,\qquad f_\ell(r) := r \frac{(\ln r)^{\ell-1}}{(\ell-1)!},\quad r>1, \quad \ell\in \N.   
\]
It was proved that 
\begin{align}\label{vol_H_1}
\vol_\ell (\mathcal{H}_\ell(r)) &\le f_\ell(r) \qquad \text{for all } ,\quad r\ge 2^\ell \quad\text{and} \quad \ell\ge 1;\\
\vol_\ell (\mathcal{H}_\ell(r)) &\ge f_\ell(r)- f_{\ell-1}(r)  \qquad \text{for all }\quad r\ge 2^\ell \quad\text{and} \quad \ell\ge 2.  \label{vol_H_2}
\end{align} 
One can easily observe that $f_\ell(r)- f_{\ell-1}(r)\ge f_\ell(r)/2$ if $r\ge e^{2(\ell-1)}$, so \eqref{vol_H_2} gives
\begin{equation}
\vol_\ell (\mathcal{H}_\ell(r)) \ge \frac{1}{2} f_\ell(r)  \qquad \text{for all }\quad r\ge e^{2(\ell-1)} .
 \label{vol_H_2a}
\end{equation}
The last inequality holds also for $\ell=1$ and $r\ge 2$. 
 and that $V^m_R = \mathcal{H}_m(R^{1/\alpha_1})$.  
Let 
\begin{equation}\nonumber
c_\ell = \alpha_1^{1-\ell}\,\frac{(\ln 2)^{\ell-1}}{(\ell-1)!}, \qquad \ell=1,\dots , m.  
\end{equation}
The above estimates give us 
\begin{equation}\label{2708-1}
\frac{1}{2} c_m R^{1/\alpha_1}(\log R)^{m-1}\le   \vol_m(V^m_R)\le c_m R^{1/\alpha_1}(\log R)^{m-1},\qquad R\ge e^{2\alpha_1(m-1)} , 
\end{equation}
since $\ell=n=m\ge 2$. 

\emph{Step 3. } Now we prove  upper and lower estimates  in the   case $n<m$. Let $\alpha_1=\ldots= \alpha_n  <\alpha_{n+1}\le \ldots \le \alpha_m$. 
We use the following relation between $\vol_{\ell+1}(V^{\ell+1}_R)$ and $\vol_{\ell}(V^{\ell}_R)$
\begin{equation}\label{ell-ell1}
\vol_{\ell+1}(V^{\ell+1}_R) = \int_1^{R^{1/\alpha_{\ell+1}}} \vol_{\ell}(V^{\ell}_{R/x^{\alpha_{\ell+1}}_{\ell+1}}) dx_{\ell+1} .
\end{equation}
First we take $\ell=n$.  The inequality  \eqref{vol_H_1}  implies 
\begin{align*}
\vol_{n+1}(V^{n+1}_R)\le &  \int_1^{R^\frac{1}{\alpha_{n+1}}} \vol_n\big(\mathcal{H}_n(R^{1/\alpha_1}/ x_{n+1}^{{\alpha_{n+1}/\alpha_1}})\big)dx_{n+1} \\
& \le \, c_n \,R^{1/\alpha_1} (\log R)^{n-1}  \int_1^\infty x^{-\alpha_{n+1}/\alpha_1}dx \,\le  \, c_n \frac{\alpha_1}{\alpha_{n+1}-\alpha_1}\,R^{1/\alpha_1} (\log R)^{n-1} ,   
\end{align*}
Iterating this argument $m-n$ times we get the  estimate 
\begin{equation}\label{VmR-up}
\vol_m(V^m_R)\le c_n \prod_{\ell=n+1}^m \frac{\alpha_1}{\alpha_{\ell+1}-\alpha_1}\;    R^{1/\alpha_1}(\log R)^{n-1}, \qquad R\ge 2^{\alpha_1 n} . 
\end{equation}
Using  \eqref{vol_H_2a} instead of  \eqref{vol_H_1} we  can prove the lower estimates in the form 
\begin{equation}\label{VmR-down}
c_n \prod_{\ell=n+1}^m \frac{\alpha_1}{\alpha_{\ell+1}-\alpha_1}\,  \le\,  \vol_m(V^m_R),  \qquad R\ge \max \{2^{\alpha_1}, e^{2\alpha_1(n-1)}\}. 
\end{equation}
This proves \eqref{VmR}. 

\emph{Step 4.} It remains to prove the estimates \eqref{red3}. By \eqref{red3a}, \eqref{VmR-up} and \eqref{VmR-down} we get 
 \begin{align*}
 \#\{k\in \N_0^m: k\in \widetilde{W}_L \} \le  2^m c_n \prod_{\ell=n+1}^m \frac{\alpha_1}{\alpha_{\ell+1}-\alpha_1} 2^{L/\alpha_1} L^{n-1}, \\
\intertext{and}
 \#\{k\in \N_0^m: k\in \widetilde{W}_L \} \ge  2^m c_n \prod_{\ell=n+1}^m \frac{\alpha_1}{\alpha_{\ell+1}-\alpha_1} (2^{\frac{1}{\alpha_1}-1}-1)\,  2^{L/\alpha_1} L^{n-1}, 
 \end{align*}
for  sufficiently large $L$. Please note that $2^{\frac{1}{\alpha_1}-1}-1>0$ since $\alpha_1< 1$. \epr

%

\begin{Prop} \label{main1seq}
	Let $\gamma= (\gamma_1 , \gamma_2, \ldots , \gamma_m)\in \N^m$, $m\in \N$,   be a multi-index such that 
	$2\le \gamma_1\le \ldots \le \gamma_m$, 
	$d= \gamma_1+ \ldots + \gamma_m$, and let $n = \max\{ i\,:\, \gamma_i=\gamma_1\}$.  
	Let $1\le p_1<p_2\le\infty$, $0<q_1,q_2\le \infty$ and  $s_1,s_2\in \R$.  If $\delta=s_1-s_2-d(\frac{1}{p_1}-\frac{1}{p_2})>0$ then 
	\begin{equation*}
	e_k\big(id: b^{\sigma_1}_{p_1,q_1}( w_\gamma)\rightarrow  b^{\sigma_2}_{p_2,q_2}( w_\gamma) \big)\, \sim\, \big( k^{-\gamma_1}(\log k)^{(n-1)(\gamma_1-1)}\big)^{\frac{1}{p_1}-\frac{1}{p_2}},
	\end{equation*}
	where $\sigma_i= s_i +\frac{m}{2}- \frac{m}{p_i}$, $i=1,2$.
\end{Prop}

\bpr
{\em Step 1.}
It is convenient to change slightly the notation. Let $0<p,q\le \infty$, $\sigma\in \R$ and let  ${\mathcal X}$ denote $\Z^m$ or $N_0$. We introduce the following sequence space,
\begin{align*}
\ell_q(2^{\nu \sigma}\ell_p({\mathcal X}, w)) &:= \left\{\right.  \lambda=\{\lambda_{\nu,\ell}\}_{\nu,\ell}: \lambda_{\nu,\ell}\in\C, \\ 
&\left. \|\lambda|l_q(2^{\nu\sigma}l_p({\mathcal X},w))\| = \bigg(\sum_{\nu=0}^\infty 2^{\nu\sigma q}\Big(\sum_{\ell\in {\mathcal X}} |\lambda_{\nu,\ell}|^p w_{\nu,\ell} \Big)^{q/p}\bigg)^{1/q}<\infty \right\} ,
\end{align*}
where  $w= (w_{\nu,\ell})_{\nu,\ell}$, $w_{\nu,\ell} >0$. As usual we write $\ell_q(\ell_p({\mathcal X}, w))$ if $\sigma=0$ and $\ell_q(2^{\nu \sigma}\ell_p({\mathcal X}))$ if $w \equiv 1$. 

By standard arguments we get 
\begin{align} \nonumber
e_k\big(id: b^{\sigma_1}_{p_1,q_1}(\ w_\gamma) &\rightarrow  b^{\sigma_2}_{p_2,q_2}( w_\gamma) \big) \\
 & \sim \,   
e_k\big( id: \ell_{q_1}(\ell_{p_1}(\Z^m))\rightarrow \ell_{q_2}(2^{\nu \sigma}\ell_{p_2}(\Z^m, w^{(\gamma)}))\big), \nonumber
\end{align}
where  $\sigma = s_2- s_1$ and  $w^{(\gamma)}= \big( w^{(\gamma)}_{\nu, n}\big)$ with  $w^{(\gamma)}_{\nu, n}= w_\gamma(Q_{\nu,n})^{1-\frac{p_2}{p_1}}$.  

{\em Step 2.} 
 Now we prove that  
\begin{align}\label{17.3-0}
e_k\Big( id: \ell_{q_1}\big(\ell_{p_1}(\Z^m)\big) & \rightarrow \ell_{q_2}\big(2^{\nu \sigma}\ell_{p_2}(\Z^m, w^{(\gamma)})\big)\Big)\, \sim  \\ 
\sim\, & e_k\Big( id: \ell_{q_1}\big(\ell_{p_1}(\N_0)\big)\rightarrow \ell_{q_2}\big(2^{-\nu \delta} \ell_{p_2}(\N_0, \widetilde{w}^{(\gamma)})\big)\Big),
\nonumber
\end{align}
where $\widetilde{w}^{(\gamma)}_\ell = \max(1,\ell\log^{1-n)}\ell)^{(\gamma_1-1)(1-\frac{p_2}{p_1})}$.
The estimate of the weight on the cubes \eqref{wQ} and the definition of the bijection $\tau$, cf. Lemma \ref{volumn}, give us 
\begin{align} \label{17.03-1}\nonumber
\Big\|\lambda \Big|  \ell_{q_2}&\big(2^{\nu \sigma}\ell_{p_2}(\Z^m, w^{(\gamma)})\big) \Big\|  = \\
& = \bigg(\sum_{\nu=0}^\infty 2^{\nu\sigma q_2}\bigg(\sum_{n\in\Z^m} |\lambda_{\nu,n}|^{p_2} w_\gamma(Q_{\nu,n})^{1-\frac{p_2}{p_1}}\bigg)^{q_2/p_2}\bigg)^{1/q_2}\sim  \nonumber\\ 
& \sim \bigg(\sum_{\nu=0}^\infty 2^{\nu(\sigma - d(\frac{1}{p_2}- \frac{1}{p_1}))q_2}\bigg(\sum_{n\in\Z^m} |\lambda_{\nu,n}|^{p_2} w_\gamma(Q_{0,n})^{1-\frac{p_2}{p_1}}\bigg)^{\frac{q_2}{p_2}}\bigg)^{\frac{1}{q_2}}\,\sim \nonumber\\
& \sim \bigg(\sum_{\nu=0}^\infty 2^{-\nu \delta q_2}\bigg(\sum_{\ell=0}^\infty |\lambda_{\nu,\tau^{-1}(\ell)}|^{p_2} w_\gamma(Q_{0,\tau^{-1}(\ell)})^{1-\frac{p_2}{p_1}}\bigg)^{\frac{q_2}{p_2}}\bigg)^{\frac{1}{q_2}}, \nonumber
\end{align}
where $\delta = s_1-s_2 - d(\frac{1}{p_1}- \frac{1}{p_2})$. 

If $L,\ell\in \N_0$ and $2^{\frac{L(d-m)}{\gamma_1-1}}L^{n-1}\le \ell <2^{\frac{(L+1)(d-m)}{\gamma_i-1}}(L+1)^{n-1}$ then 
\begin{equation}\label{red2a}
w_\gamma (Q_{0,\tau^{-1}(\ell)}) \,\sim\,  2^{L(d-m)}\, \sim \, (\ell\log^{1-n}\ell)^{\gamma_1-1}\, ,
\end{equation}
cf. \eqref{red1}. By what we have already proved 
\begin{align}\label{17.3-2}\nonumber
\Big\|\lambda \Big|\ell_{q_2}\big(2^{\nu \sigma}&\ell_{p_2}(\Z^m, w^{(\gamma)})\big) \Big\|   \, \sim \, \\
&  \sim \bigg(\sum_{\nu=0}^\infty 2^{-\nu \delta q_2}\bigg(\sum_{\ell=0}^\infty |\lambda_{\nu,\tau^{-1}(\ell)}|^{p_2} 
\max(1,\ell\log^{1-n}\ell)^{(\gamma_1-1)(1-\frac{p_2}{p_1})}\bigg)^{\frac{q_2}{p_2}}\bigg)^{\frac{1}{q_2}}. \nonumber
\end{align}
This justifies the equivalence \eqref{17.3-0} since the estimate 
\[
\Big\|\lambda \Big|\ell_{q_1}\big(\ell_{p_1}(\Z^m)\big) \Big\|   \, \sim \, \bigg(\sum_{\nu=0}^\infty \bigg(\sum_{\ell=0}^\infty |\lambda_{\nu,\tau^{-1}(\ell)}|^{p_1} \bigg)^{\frac{q_1}{p_1}}\bigg)^{\frac{1}{q_1}}  
\]
is obvious.

{\em Step 3.}  We prove the upper estimate of the entropy numbers 
$$
e_k\Big( \id: \ell_{q_1}\big(\ell_{p_1}(\N_0)\big)\rightarrow \ell_{q_2}\big(2^{-\nu \delta} \ell_{p_2}(\N_0, \widetilde{w}^{(\gamma)})\big)\Big)\, 
\le \, \big( k^{-\gamma_1}(\log k)^{(n-1)(\gamma_1-1)}\big)^{\frac{1}{p_1}-\frac{1}{p_2}}.
$$
Let us consider the projection 
$P_\nu: \ell_{q_1}\big(\ell_{p_1}(\N_0)\big) \rightarrow \ell_{p_1}(\N_0)$
onto the $\nu$-th vector-coordinate,  and  the embedding operator $E_\nu: \ell_{p_2}(\N_0) \rightarrow \ell_{q_2}\big(2^{-\nu \delta} \ell_{p_2}(\N_0)\big)$, 
\[
\big(E_\nu(y)\big)_{\mu, \ell}= \begin{cases}
y_\ell & \text{if}\qquad  \mu = \nu\\
0 & \text{otherwise}, 
\end{cases}
\qquad y \in  \ell_{p_2}\big(\N_0, \big)
\]
It is  obvious that 
\begin{equation}\nonumber
\|P_\nu\| = 1 \qquad \text{and}\qquad \|E_\nu\| = 2^{-\nu\delta}\, .
\end{equation}
Let $D_\gamma$ denote the diagonal operator $D_\gamma: \ell_{p_1}(\N_0) \rightarrow \ell_{p_2}(\N_0)$ generated by the sequence
$\sigma_\ell =  (\widetilde{w}_\ell^{(\gamma)})^{\frac{1}{p_2}}$ i.e. 
$$
(D_\gamma(\lambda))_\ell = (\widetilde{w}^{(\gamma)}_\ell)^{\frac{1}{p_2}} \lambda_\ell\, .
$$ 
Then
$$
\id = \sum_{\nu=0}^\infty \id_\nu, \qquad \text{where} \qquad
\id_\nu=E_\nu D_\gamma P_\nu \,.
$$
The multiplicativity of the  entropy numbers yields
$$
e_k(\id_\nu) \leq c \,2^{-\nu\delta} e_k(D_\gamma).
$$
Now using  Theorem \ref{Kuhn} 
we get 
 $$e_k(D_\gamma)\sim k^{-(1/p_1-1/p_2)}(k\log^{1-n} k)^{-(\gamma_1 - 1)(1/p_1-1/p_2)}.$$
So
$$
e_k(\id_\nu) \leq c \,2^{-\nu\delta} \big(k^{-\gamma_1 }(\log k)^{(n-1)(\gamma_1 - 1)}\big)^{1/p_1-1/p_2}
$$
with a constant $c$ independent of $\nu$ and $k$. Now taking $\omega_k=(k^{\gamma_1 }(\log k)^{(1-n)(\gamma_1 - 1)})^{1/p_1-1/p_2}$ we have
$$
L_{\omega}^{(e)}(\id_\nu)=\sup _k \, \omega_k e_k(\id_\nu) \leq c
\,2^{-\nu\delta}.  
$$
Since quasi-norm $L_{\omega}^{(e)}$ is equivalent to an $r$-norm for some $r$, $0<r \leq
1$, we arrive at 
$$
L_{\omega}^{(e)}(\id)^r \leq c \sum_{\nu=0}^\infty
L_{\omega}^{(e)}(\id_\nu)^r \leq c \sum_{\nu=0}^\infty 2^{-\nu \delta r} < \infty \, , 
$$
which proves the upper estimate.

{\em Step 4.}  In this step we  estimate the entropy numbers from below. For any given $k \in \N$ we consider $k$-dimensional vector spaces $\ell^k_{p_i}$ and the following commutative diagram
\begin{equation}\nonumber
\begin{CD}
\ell_{p_1}^{k}@>{T}>>\ell_{q_1}\big(\ell_{p_1}(\N_0)\big)\\
@V{\Id}VV@VV{\id}V\\
\ell_{p_2}^{k} @<{S}<<\ell_{q_2}\big(2^{-\nu\delta}\ell_{p_2}(\N_0, \tilde{w}^{(\gamma)})\big).\,
\end{CD}
\end{equation}
Here the operators $S$ and $T$ are defined by  
\begin{align*}
\big(T(\xi_1, ... , \xi_k)\big)_{\nu,\ell} = 
\begin{cases} 
\xi_{\ell + 1-k}  & \text{if}\;  \nu=0 \;\text{and}\; k \leq l \leq 2k-1,\\
0  &\text{otherwise} 
\end{cases}
\end{align*}
and
$$
S((\lambda_{\nu,\ell})_{\nu,\ell} ) = (\lambda_{0,k}, \cdots ,\lambda_{0,2k-1}). 
$$

The norms of the above operators have the obvious estimates 
$$
\|T\| \leq 1    \quad {\rm and } \quad   \|S\| \leq (\widetilde{w}^{(\gamma)}_k)^{-\frac{1}{p_2}} \, .
$$
Using Sch\"utt's description of asymptotic behaviour of entropy numbers for embeddings between the finite dimensional 
spaces $ \ell_p^k$ - see \cite{Sch}, and again the multiplicativity of the entropy numbers
we get,   
\begin{align*}
c k^{-(\frac{1}{p_1}-\frac{1}{p_2})} & \leq  e_k(\id: \ell_{p_1}^k \rightarrow \ell_{p_2}^k) \nonumber \\[0.1cm]
& \leq \|S\|  \,e_k(\id: \ell_{q_1}(\ell_{p_1}(\N_0)) \rightarrow      \nonumber
\ell_{q_2}(2^{-\nu \delta}\ell_{p_2}(\N_0, \widetilde{w}^{(\gamma)}))) \,\|T\| \le \\
& (k\log^{1-n} k)^{(\gamma_1-1)(\frac{1}{p_1}-\frac{1}{p_2})}\,e_k(\id: \ell_{q_1}(\ell_{p_1}(\N_0)) \rightarrow      \nonumber
\ell_{q_2}(2^{-\nu \delta}\ell_{p_2}(\N_0, \widetilde{w}^{(\gamma)})))\,  
\end{align*}
with some constant $c$ independent of $k$. This proves the proposition. \epr


\begin{Prop} \label{main1function}
	Let $\gamma= (\gamma_1 , \gamma_2, \ldots , \gamma_m)\in \N^m$, $m\in \N$,   $2\le \gamma_1\le \ldots \le \gamma_m$, 
	$d= \gamma_1+ \ldots + \gamma_m$, and let $n = \max\{ i\,:\, \gamma_i=\gamma_1\}$. 	
	Let $1\le p_1<p_2\le\infty$, $0<q_1,q_2\le \infty$ and  $s_1,s_2\in \R$.  If $s_1-s_2-d(\frac{1}{p_1}-\frac{1}{p_2})>0$ then 
	\begin{equation*}
	e_k\Big(\id :R_G B^{s_1}_{p_1,q_1}(\R^m, w_\gamma)\rightarrow R_G B^{s_2}_{p_2,q_2}(\R^m, w_\gamma)\Big)\, \sim\,  \big(k^{-\gamma_1}(\log k)^{(n-1)(\gamma_1-1)}\big)^{\frac{1}{p_1}-\frac{1}{p_2}}.
	\end{equation*}
\end{Prop}

\bpr
Proposition \ref{main1seq} and the wavelet decomposition of the spaces give us 
\begin{equation}\label{bweit}
e_k\Big(\Id :B^{s_1}_{p_1,q_1}(\R^m, w_\gamma)\rightarrow B^{s_2}_{p_2,q_2}(\R^m, w_\gamma)\Big)\, \sim\,  \big(k^{-\gamma_1}(\log k)^{(n-1)(\gamma_1-1)}\big)^{\frac{1}{p_1}-\frac{1}{p_2}}.
\end{equation}
So it remains to show that we have the same estimate for the $G(\gamma)$-invariant subspaces. The operator 
$$
P f (x)= \frac{1}{|G(\gamma)|}\sum_{g\in G(\gamma)} f(g(x)) 
$$
is a bounded projection of $B^{s}_{p,q}(\R^m, w_\gamma)$ onto $R_G B^{s}_{p,q}(\R^m, w_\gamma)$ since $G(\gamma)$ is a finite group of linear  isometries.  So 
the inequality 
\[
e_k\Big(\id:R_G B^{s_1}_{p_1,q_1}(\R^m, w_\gamma)\rightarrow R_G B^{s_2}_{p_2,q_2}(\R^m, w_\gamma)\Big)\, \le\, C \big(k^{-\gamma_1}(\log k)^{(n-1)(\gamma_1-1)}\big)^{\frac{1}{p_1}-\frac{1}{p_2}}
\]
follows from \eqref{bweit} and  the following  commutative diagram 
\begin{equation}\nonumber
\begin{CD}
R_GB^{s_1}_{p_1,q_1}(\R^m, w_\gamma) @ >{\id}>> R_GB^{s_2}_{p_2,q_2}(\R^m, w_\gamma)\\
@V{id}VV                                          @AA{P}A \\
B^{s_1}_{p_1,q_1}(\R^m, w_\gamma) @ >>{\Id}> B^{s_2}_{p_2,q_2}(\R^m, w_\gamma). 
\end{CD}
\end{equation}

To prove the opposite  inequality we use the wavelet decomposition, cf. Appendix B. The group $G(\gamma)$ divides $\R^m$ into the finite sum of cones, that have pairwise disjoint interiors. Let us choose one of those cones and denote it by $\widetilde{\mathcal{C}}$. 
Moreover let $\eta$ be a function belonging to $C^\ell(\Rd)$, $\ell>s_1$, such that $\supp \eta \subset \widetilde{\mathcal{C}}$ and $\eta(x)=1$ if $x\in {\mathcal{C}} = \{x\in \widetilde{\mathcal{C}}: \dist (x, \partial\widetilde{\mathcal{C}})>\varepsilon\}$, for some fixed sufficiently small $\varepsilon>0$. 
We consider the family  $\mathcal{K} = \{k\in \Z^m:\; \supp \phi_{0,k} \subset \mathcal{C}\}$. The group $G(\gamma)$ is finite therefore for any $0<c_1<c_2$ we can find   $L_0$ such that for any $L\ge L_0$, $L\in \N$, we have 
\begin{align}\label{below1}
\#\{k\in \mathcal{K}: & \;c_1 2^{L(d-m)}\le w(Q_{0,k})\le c_2 2^{L(d-m)} \} \, \sim\,  \\ 
\sim\, &  \#\{k\in \Z^d:   c_1 2^{L(d-m)}\le w(Q_{0,k})\le c_2 2^{L(d-m)} \} .  \nonumber
\end{align}   
Please note that $L$ large means that the cubes $Q_{0,k}$ are located far from the origin, cf. \eqref{wQ}. The above considerations and Lemma \ref{volumn} yield  the existence of the bijection $\sigma: \mathcal{K} \rightarrow \N_0$ such that 
\begin{equation}\label{red1bis}\nonumber
w_\gamma (Q_{0,k}) \,\sim\,  2^{L(d-m)} \qquad \Longleftrightarrow\qquad  \sigma(k)\, \sim \, 2^{\frac{L(d-m)}{\gamma_1-1}}L^{n-1} \, .
\end{equation}

Let $v_{\gamma}(\ell) =  w_\gamma (Q_{0,\sigma^{-1}(\ell)})$.  Then \eqref{red2a} and \eqref{below1} imply
\begin{equation}\label{below1a}
v_\gamma(\ell)  
\, \sim \, (\ell\log^{1-n}\ell)^{\gamma_1-1}\, .
\end{equation}  

For  further arguments we need three linear bounded  operators: $T:\ell_{p_1}(\N_0,v_\gamma)\rightarrow \ell_{q_1}(\ell_{p_1}(\Z^m, w_\gamma))$, $M_\eta: B^{s_2}_{p_2,q_2}(\R^m, w_\gamma)\rightarrow B^{s_2}_{p_2,q_2}(\R^m, w_\gamma)$ and $S:\ell_{q_2}(\ell_{p_2}(\Z^m, w_\gamma))\rightarrow \ell_{p_2}(\N_0,v_\gamma)$. The operators are defined in the following way:
\begin{align*}
(T\lambda)_{j,k} &\, = \begin{cases}
\lambda_{\sigma(k)}& \quad\text{if}\quad k\in \mathcal{K}\quad \text{and}\quad j=0,\\
0 & \quad \text{otherwise},
\end{cases}
\qquad \lambda\in \ell_{p_1}(\N_0,v_\gamma);\\
M_\eta(f) &\, = |G(\gamma)| \eta \cdot f, \qquad f\in  B^{s_2}_{p_2,q_2}(\R^m, w_\gamma) ; \\
S(\lambda)_\ell &\, = \lambda_{0, \sigma^{-1}(\ell)} \qquad \text{it}\qquad \lambda\in  \ell_{q_2}(\ell_{p_2}(\Z^m, w_\gamma))\, .
\end{align*}

Using these operators we can construct the following commutative diagram 
\begin{equation}\label{diag230317}\nonumber 
\begin{CD}
\ell_{p_1}(\N_0,v_\gamma)@>{T}>> \ell_{q_1}(\ell_{p_1}(\Z^m, w_\gamma))     @>{\mathcal{W}^{-1}}>> B^{s_1}_{p_1,q_1}(\R^m, w_\gamma) @ >{P}>>  R_GB^{s_1}_{p_1,q_1}(\R^m, w_\gamma) \\ 
@V{\Id}VV     &{} &{} &{}  &{}                                                        @VV{\id}V\\
\ell_{p_2}(\N_0,v_\gamma)@<S<<  \ell_{q_2}(\ell_{p_2}(\Z^m, w_\gamma)) @<{\mathcal{W}}<< B^{s_2}_{p_2,q_2}(\R^m, w_\gamma) @ <{M}<<  R_GB^{s_2}_{p_2,q_2}(\R^m, w_\gamma).
\end{CD}
\end{equation}

Here   $\mathcal{W}$ is the isomorphism defined by the wavelet basis, cf. Appendix B. It follows from the above diagram that  
\begin{equation}\label{below2}
e_k(\Id: \ell_{p_1}(\N_0,v_\gamma)\rightarrow \ell_{p_2}(\N_0, v_\gamma))\le C  e_k (\id).
\end{equation} 
But 
\begin{equation}\label{below3}
e_k(\Id: \ell_{p_1}(\N_0,v_\gamma)\rightarrow \ell_{p_2}(\N_0,v_\gamma))\, \sim \, 
e_k(D_\gamma: \ell_{p_1}(\N_0)\rightarrow \ell_{p_2}(\N_0)),
\end{equation} 
where  $D_\gamma$ denote the diagonal operator  generated by the sequence
$\sigma_\ell =  (\ell\log^{1-n}\ell)^{(\gamma_1 - 1)(1/p_2-1/p_1)}$ i.e. 
$$
(D_\gamma(\lambda))_\ell = (\ell\log^{1-n}\ell)^{(\gamma_1 - 1)(1/p_2-1/p_1)} \lambda_\ell\, .
$$ 
Using once more  K\"uhn's results from \cite{Kuh} or \cite{Kuh2} and \eqref{below1a}-\eqref{below3} we have
 $$ k^{1/p_1-1/p_2}(k\log^{1-n} k)^{(\gamma_1 - 1)(1/p_1-1/p_2)} \sim e_k(D_\gamma)\le  e_k (\id).$$
This proves the proposition.
\epr

Summarizing, we have the following theorem.
\begin{T}\label{mainT1}
			Let $\gamma= (\gamma_1 , \gamma_2, \ldots , \gamma_m)\in \N^m$, $m\in \N$,   be a multi-index such that $2\le \gamma_1\le \ldots \le \gamma_m$, 
	$d= \gamma_1+ \ldots + \gamma_m$, and let $n = \max\{ i\,:\, \gamma_i=\gamma_1\}$.	

	Let $1 < p_1<p_2 \leq \infty$, $0<q_1,q_2\le \infty$ and  $s_1,s_2\in \R$.  If  $\delta=s_1-s_2-d(\frac{1}{p_1}-\frac{1}{p_2})>0$ then 
	\begin{align}\label{mainT1B}
	e_k\Big( \id: R_{\gamma} B^{s_1}_{p_1,q_1}(\Rd)\rightarrow R_{\gamma} B^{s_2}_{p_2,q_2}(\Rd)\Big)\, \sim\, 
	\big(k^{-\gamma_1}(\log k&)^{(n-1)(\gamma_1-1)}\big)^{\frac{1}{p_1}-\frac{1}{p_2}} \\
	\intertext{and}
	e_k\Big( \id: R_{\gamma} H^{s_1}_{p_1}(\Rd)\rightarrow R_{\gamma} H^{s_2}_{p_2}(\Rd)\Big)\, \sim\, 
	\big(k^{-\gamma_1}(\log k&)^{(n-1)(\gamma_1-1)}\big)^{\frac{1}{p_1}-\frac{1}{p_2}}.  \label{mainT1H}
	\end{align}
\end{T}
\begin{Rem}{\rm 
		If $m=1$ then the spaces  $R_{\gamma} B^{s}_{p_,q}(\Rd)$ and $R_{\gamma} H^{s}_{p}(\Rd)$ consists  of radial distributions and the estimates \eqref{mainT1B}-\eqref{mainT1H} coincides with the estimates for radial functions proved in \cite{KLSS}. If $\gamma_i\not= \gamma_j$ for any $i\not=j$ then 
		\[
		e_k(\id) \, \sim \, k^{-\gamma_1(\frac{1}{p_1}-\frac{1}{p_2})} .
		\]
		So the asymptotic behaviour  is the same as for the corresponding radial subspaces defined on the block with the lowest dimension. 
		On the other hand if $\gamma_1=\ldots =\gamma_m$ then 
	            \[
		e_k(\id) \, \sim \,  \big(k^{-\gamma_1}(\log k)^{(m-1)(\gamma_1-1)}\big)^{\frac{1}{p_1}-\frac{1}{p_2}} .\]	 
		In that case the sequence of the entropy numbers for spaces of block-radial functions goes asymptotically to zero  slower  than the entropy numbers for radial function defined on any block.  
	}\end{Rem}

	
		\section{Block-radial bounded states of Schr\"odinger type operators}
		\label{neg_spec}
		
		
		The interest in studying the 'negative' spectrum (bound states) comes from quantum
		mechanics, generalizing the classical hydrogen operator,
		\beq\label{hydro}\nonumber
		H=-\Delta-\frac{c}{|x|}\;,\quad c>0\;,
		\eeq
		in $L_2(\R^3)$. Thus 'potentials' $V(x)$ with 
		$V(x)\sim |x|^{-a}$, $a>0$, are of peculiar interest. These potentials have not only   local singularities and some decay properties at  infinity but also they are radial. Here, more generally, we want to consider the 'potentials' which have block-radial symmetry. We want to estimate the number of negative eigenvalues  that corresponds to block-radial eigenvectors, but  first we briefly describe the general setting.  
		
		\subsection{The Birman-Schwinger principle and the Carl inequality}
		We adapt the Birman-Schwinger principle as described in \cite{Schech} and
		\cite{Simon} to our concrete situation. Let $A$ be a self-adjoint
		positive-definite operator and let $B$ be a symmetric relatively compact operator in
		the Hilbert space $\cal H$. Let $\psigma$ denote the point spectrum and $\esigma$ denote
		the essential spectrum of a self-adjoint operator. Then the
		eigenvalues $\{\mu_k\}_k$ of $BA^{-1}$ are real, and $(B A^{-1})^\ast =
		A^{-1}B$ is the adjoint operator after extension by continuity from $\dom(B)$
		to $\cal H$. Furthermore, the operator $A+B$ with $\dom(A+B) = \dom(A)$, is
		self-adjoint, with $\esigma(A+B) = \esigma(A)$, and
		\begin{align}
		\# \left\{\psigma(A+B) \cap (-\infty,0]\right\} & = \# \left\{\sigma (A+B) \cap
		(-\infty,0]\right\}\nonumber\\
		&= \# \left\{k\in\N : \mu_k(B A^{-1}) \leq -1 \right\} < \infty.
		\label{6.80}
		\end{align}
		This is usually called the {\em Birman-Schwinger principle}. It goes back to
		\cite{birman,schwinger}, proofs may be found in
		\cite[Chapter~7]{Simon} and \cite[Chapter~8, \S~5]{Schech}. A short
		description has also been given in \cite[Section~5.2.1, p.~186]{ET}. Our
		formulation is different and adapted to our later needs. 
		
		Our approach is based upon  the relation between  the eigenvalues of a compact operator and its entropy numbers described by  the  Carl inequality. If $\big(\lambda_k(T)\big)_{k\in \N}$ is a decreasing sequence of all non-zero eigenvalues of  a compact operator $T$, repeated according to their algebraic multiplicities then the following inequality
		\begin{equation}\label{Carl}
		|\lambda_k(T)|\;\le\; \sqrt{2}\, e_k(T)\,
		\end{equation}
		holds, cf. \cite{Carl, Ca-Trie},  \cite[Theorem 1.3.4]{ET}.
		Using \eqref{Carl} with $T=BA^{-1}$ one
		obtains by \eqref{6.80} 
		\begin{equation*}
		\#\left\{\psigma(A+B)\cap (-\infty,0]\right\}\leq \# \left\{k\in\N: \sqrt{2}
		e_k\left(BA^{-1}\right)\geq 1\right\}.
		\end{equation*}
		This entropy version of the Birman-Schwinger principle appeared first in
		\cite[Theorem~2.4]{HT1}, cf. also \cite[Corollary, p.~186]{ET}.

		We shall concentrate on the special case when $B= -V$  is a multiplication
		operator where (in a slight abuse of notation) $V$ is a 
		nonnegative measurable and $SO(\gamma)$-invariant function, finite a.e., typically belonging to some
		space $L_r(\R^d)$. 
		
		We turn to study the behaviour of the part of negative spectrum of the self-adjoint unbounded operator 
		\begin{equation}\label{Hgamma}
		H_{s,\theta,\beta} = (\theta\Id -\Delta)^{s/2} - \beta V\qquad\mbox{as}\quad
		\beta\rightarrow\infty ;\qquad 0<\theta\le 1,  
		\end{equation}
		corresponding to the $SO(\gamma)$-invariant eigenfunctions. We assume that $s>0$, $\beta>0$ and $V\ge 0$ is an  $SO(\gamma)$-invariant potential. The operator $H_{s,\theta,\beta}$ is a bounded below, self-adjoint operator in
		$L_2({\R}^d)$ with the domain $D(H_{s,\theta,\beta})\;=\; H^s_2({\R}^d)$.
	   Let $\psigmag$ denote the part of the point spectrum of the operator   
		\eqref{Hgamma} that corresponds to the $SO(\gamma)$-invariant eigenfunctions.  By the Birman-Schwinger principle with ${\cal H}= R_{\gamma}L_2$ as  the basic space we get 
		\beq\nonumber
		\# \{\psigmag(H_{s,\theta,\beta}) \cap (-\infty, 0]\}\leq \# \left\{k\in\N:
		\sqrt{2\,}\;e_k\left(V^\frac12 (\theta\Id -\Delta)^{-s/2} V^\frac12\right)
		\geq \lambda^{-1}\right\}.\label{nesp-1}
		\eeq
		
		Thus we should consider the compactness and asymptotic behaviour of entropy numbers of the operators  
		$ V_2\Delta_\theta^{-s/2}V_1$ where $V_1$, $V_2$ are  positive block-radial functions and $\Delta_\theta:=\theta\Id -\Delta$.

		\begin{Lem}  \label{oper4} 
				Let $\gamma= (\gamma_1 , \gamma_2, \ldots , \gamma_m)\in \N^m$, $m\in \N$,   be a multi-index such that $2\le \gamma_1\le \ldots \le \gamma_m$, 
	$d= \gamma_1+ \ldots + \gamma_m$, and let $n = \max\{ i\,:\, \gamma_i=\gamma_1\}$.
			Let $1 \le p \le  \infty$,  $s> 0$,    $V_1\in R_{\gamma}L_{r_1}(\R^d)$, $V_2\in R_{\gamma}L_{r_2}(\R^d)$ and 
			$\frac{s}{d}> \frac{1}{r_1} + \frac{1}{r_2}>0$. If $p' < r_1 \le  \infty$ and $p \le  r_2 \le \infty$ 
			then the operator 
			$$
			V_2(\Delta_\theta)^{-s/2}V_1\, :\, R_{\gamma}L_p({\R}^d) \rightarrow R_{\gamma}L_p({\R}^d)
			$$ 
			is compact.  
			Moreover its eigenvalues and entropy numbers satisfy the following estimate 
			\begin{align}\label{op-ent1a}
			\lambda_k\big(V_2\Delta_\theta^{-s/2}V_1\big)\;\le &\;   \sqrt{2}\, e_k\big(V_2\Delta_\theta^{-s/2}V_1\big)\;\le \\ & \le  C \big(k^{-\gamma_1} (\log k)^{(n-1)(\gamma_1 -1) (\frac{1}{r_1}+\frac{1}{r_2})} \|V_1|L_{r_1}(\R^d)\| \, \|V_2|L_{r_2}(\R^d)\| \, .\nonumber
			\end{align}
		\end{Lem}
		
		\bpr
		The reasoning goes by standard factorization
		\begin{equation}\label{fact}
		\begin{tikzcd}
		R_{\gamma}L_p({\R}^d) \arrow[rr, "V_2\Delta_\theta^{-s/2}V_1"]   \arrow[d, "V_1"] &    & R_{\gamma}L_p({\R}^d)  \\
		R_{\gamma}B^0_{r,\infty}({\R}^d)  \arrow[r, "\Delta_\theta^{-s/2}"]  &  
		R_{\gamma}B^s_{r,\infty}({\R}^d) \arrow[r, "\Id"]  
		& R_{\gamma}B^0_{t,1}({\R}^d),  \arrow[u, "V_2"]\\
		\end{tikzcd}
		\end{equation} 
		where $\frac{ 1}{ r} = \frac 1 p + \frac{ 1}{ r_1}$, $\frac{ 1}{ t}=\frac{ 1}{ p} - \frac{ 1}{ r_2}$ and  $V_i$ denotes  the operator of multiplication by function $V_i$, $i=1,2$.
		The operator  $(\Delta_\theta)^{-s/2}:R_\gamma B^0_{r,\infty}({\R}^d)\rightarrow R_\gamma B^s_{r,\infty}({\R}^d)$ is an isomorphism and  the Sobolev embedding $\Id :R_\gamma B^s_{r,\infty}({\R}^d)\rightarrow R_\gamma B^0_{t,1}({\R}^d)$ is compact since  $s >\frac{ d}{ r_1}  + \frac{ d}{ r_2} = \frac{ d}{ r}  -\frac{ d}{ t}>0$.
		Moreover,  $\frac 1 p=\frac 1 r - \frac{1}{r_1}=\frac 1 t + \frac{1}{r_2}$ therefore  by H\"older inequality and elementary embeddings $B^0_{\rho,1}(\R^d) \hookrightarrow L_\rho (\R^d)  \hookrightarrow B^0_{\rho,\infty}$ the operators    $V_1:R_\gamma L_p({\R}^d)\rightarrow R_{\gamma}B^0_{r,\infty}({\R}^d)$,  $V_2:R_{\gamma}B^0_{t,1}({\R}^d)\rightarrow R_\gamma L_{p}({\R}^d)$ are bounded.  The lemma follows from Theorem \ref{mainT1} and the Carl inequality \eqref{Carl}.  
		\epr
		
		\br{op-rem1}{
			The constant $C$ in 
			\eqref{op-ent1a} depends on $s$, $p$, $r_1$, $r_2$ and $\theta$. It follows  from the proof that $C=c \| (\Delta_\theta)^{-s/2}\|$ , where $c$ is independent of $\theta$.
			
			In the similar way one can estimate eigenvalues and entropy numbers of the operators  $V\Delta_\theta^{-s/2}$ and $\Delta_\theta^{-s/2}V$.}
		\er

		\subsection {The negative  spectrum of Schr\"odinger type operators}

		We are interested  in a number of negative eigenvalues of $H_{s,\beta,\theta}$ with $SO(\gamma)$-invariant eigenfunctions. We put
		\begin{equation}\label{Nbeta}\nonumber
		N_{\gamma, \beta} =\# \Big\{ \lambda\, :\; \lambda\le 0\, ,\quad  H_{s,\theta,\beta}f=\lambda f\, ,\quad f \in R_{\gamma}L_2({\R}^d)\, ,\; f\not= 0\,  \Big\}\, .
		\end{equation}
		
		The main theorem of this section reads as follows
		\bt\label{op-main} 
			Let $\gamma= (\gamma_1 , \gamma_2, \ldots , \gamma_m)\in \N^m$, $m\in \N$,   be a multi-index such that $2\le \gamma_1\le \ldots \le \gamma_m$, 
	$d= \gamma_1+ \ldots + \gamma_m$, and let $n = \max\{ i\,:\, \gamma_i=\gamma_1\}$.
		Let $V$ be a nonnegative block-radial function invariant with respect to $SO(\gamma)$ such that   $\|V|L_r({\R}^d)\|=1$,  $1 < r<\infty$.
		Let $0< \theta< 1$, $s>0$   and $\beta>0$. We assume that  $\frac s d > \frac 1 r$.
		If $H_{s,\theta,\beta}$ is the  operator defined by \eqref{Hgamma} with domain
		$H^s_2({\R}^d)$  then
		\[
		N_{\gamma,\beta}\;\le\; c\, \beta^{\frac{r}{ \gamma_1}} \big(\log \beta\big)^{(n-1)\frac{\gamma_1-1}{\gamma_1}} .
		\]
		Moreover, if  there exist  $\delta>0$, $\varepsilon>0$ and  $1<\rho_1, \dots , \rho_m$   
		such that
		\begin{equation}\label{loweas}\nonumber
		V(x)\ge \varepsilon \qquad \text{if} \qquad x\in A= \{x\in \R^d:\; \rho_i\le r_i(x)\le \rho_i+\delta ,\; i=1,\ldots , m \}
		\end{equation}
		then
		\begin{equation}
		\label{eig4}\nonumber
		c \beta^{\frac m s}\le  N_{\gamma,\beta}\ .
		\end{equation}
		\et
		
		\bpr {\em Step 1.} The number $\lambda$ is  a negative eigenvalue of $H_{s,\theta,\beta}$ with block-radial eigenfunction if and only if $\lambda$ belongs to the spectrum of $H_{s,\theta,\beta}$ regarded as an operator in $R_{\gamma}L_2({\R}^d)$ with domain $R_{\gamma}H^s_2({\R}^d)$.  Thus it is sufficient to consider  the last operator.
		The operator $\beta \sqrt{V}\Delta^{-s/2}_\theta\sqrt{V}$ is compact in $R_{\gamma}L_2({\R}^d)$.  This follows from the factorization \eqref{fact} with $p=2$ and $r_1=r_2=2r$. 
		By the Birman-Schwinger principle the operator $H_{s,\beta,\theta}$ is self-adjoint with the same domain as $\Delta^{s/2}_\theta$, cf.   \cite[Proposition 5.4.1]{ET}. Moreover
		\[
		N_{\gamma,\beta}\;\le\;\#\big\{k\in \N:\; \sqrt{2}e_k\big(\beta \sqrt{V}\Delta_\theta^{-s/2}\sqrt{V}\big)\ge 1\big\}\, ,
		\]
		So the upper estimate  follows from Lemma \ref{oper4}, since
		\[
		e_k\big(\beta \sqrt{V} \Delta^{-s/2}_\theta \sqrt{V}  \big)\;\le\; C\;
		 \big(k^{-\gamma_1} (\log k)^{(n-1)(\gamma_1 -1)} \big)^{\frac{1}{r}} \beta\, ,
		 \]
		cf. \reff{op-ent1a}. 
		
		{\em Step 2.} Now we prove the estimate from below. We use  the Max-Min principle and the method of atomic decompositions, cf. Appendix A. We have
		\begin{eqnarray*}
			\big(H_{s,\beta,\theta}f,f)&\sim& \|f|H^{s/2}_2({\R}^d)\|^2 -\beta (Vf,f)_{L_2({\R}^d)} \\
			&\sim&\|f|B^{s/2}_{2,2}({\R}^d)\|^2 -\beta \|f\sqrt{V}|L_2({\R}^d)\|^2, \qquad f\in H^{s/2}_2(\R^d)\, . \nonumber
		\end{eqnarray*}

		{\em Substep 2.1.}
		Let $\eta \in C^\infty_0(\R)$ be a smooth function such that $\supp \eta \subset [0,\delta]$, $0\le \eta(t)\le 1$ and $\eta(t)=1$ if $t\in [\frac \delta 4, \frac{3\delta}{4}]$. Let $\eta_j(t)= 
		\eta(2^jt)$, $j=0,1,2,\dots $. We put 
		\begin{align}\nonumber
		r^{(i)}_{j,\nu} = \rho_i + \nu 2^{-j}\delta ,\qquad 
		\widetilde{r}^{(i)}_{j,\nu} = r^{(i)}_{j,\nu}+2^{-j-2}\delta\quad \text{and}\quad \bar{r}^{(i)}_{j,\nu} = r^{(i)}_{j,\nu}+3\cdot 2^{-j-2}\delta,  
		\end{align}
		where $ \nu =0,1,2,\ldots 2^j -1 $. 
		
		We choose the following functions 
		\begin{align}\nonumber
		\psi_{j,\widetilde{\nu}}(x)=  2^{-j\frac{s- d}{2}} \prod_{i=1}^m\eta^{(i)}_{j,\nu_i}(r_i(x)), \;  x\in \R^d\, \quad\text{where}\quad \eta^{(i)}_{j,\nu_i}(t) = \eta_j(t-r^{(i)}_{j,\nu_i}) ,
		\intertext{and}
		j\in \N,\quad  \widetilde{\nu}=(\nu_1,\ldots ,\nu_m),\quad   \nu_i=0,\ldots , 2^j -1, \quad i=1,\ldots ,m . \nonumber
		\end{align}
		Any function $\psi_{j,\widetilde{\nu}}$ is  smooth and  $SO(\gamma)$-invariant.  Such a function is   supported in the set  
		$$ 
		\{x\in \R^d:\; r^{(i)}_{j,\nu_i} \le r_i(x)\le r^{(i)}_{j,\nu_i+1} ,\; i=1,\ldots , m \} \subset A 
		$$ 
		and it  takes value  $2^{-j\frac{s- d}{2}}$ on the set  $A_{j,\widetilde{\nu}}=\{x\in \R^d:\; \widetilde{r}^{(i)}_{j,\nu_i} \le r_i(x)\le \bar{r}^{(i)}_{j,\nu_i} ,\; i=1,\ldots , m \} $. 
		Moreover there is a positive constant $C$ such that
		\begin{equation}\label{op-below2}
		\Big| \partial^\alpha \psi_{j,\widetilde{\nu}}(x)\Big|\;\le\; C\, 2^{-j(\frac{s}{2}-|\alpha|-\frac d 2)},\qquad |\alpha|\le s+1 .
		\end{equation}

		Let $\{x^{(j)}_{\tilde{k},\tilde{\ell}}\}$ be $(2^{-j}\delta)$-discretization described in the appendix.  Then  
		the balls $\{B(x^{(j)}_{\tilde{k},\tilde{\ell}}, 2\delta\, 2^{-j})\}_{\tilde{k},\tilde{\ell}}\}$ form a uniformly locally finite covering of 
		$\R^d$, cf. Remark \ref{wwwr} ibidem. So there exists a resolution of unity 
		$\big( \varphi_{j,\tilde{k},\tilde{\ell}}\big)_{\tilde{k},\tilde{\ell}}$  related to this covering  such that
		\begin{equation}\label{op-below2a}
		\Big| \partial^\alpha\varphi_{j, \tilde{k},\tilde{\ell}}(x)\Big|\;\le\; C\, 2^{j|\alpha|}, \qquad |\alpha|\leq s+1\ .
		\end{equation}
		It follows from \eqref{op-below2} and \eqref{op-below2a}  that functions
		$$
		a_{j,\tilde{k},\tilde{\ell}}(x)=\varphi_{j,\tilde{k},\tilde{\ell}}(x)\psi_{j,\widetilde{\nu}}(x)
		$$
		are $(\frac{s}{2},2)$-atoms and
		$$
		\psi_{j,\widetilde{\nu}}\;=\; \sum_k\sum_\ell a_{j,\tilde{k},\tilde{\ell}}\,
		$$
		is the atomic decomposition of  $\psi_{j,\widetilde{\nu}}$.
		
		Let $k_{i,\widetilde{\nu}}$ be an integer such that $2^{-j}(k_{i,\widetilde{\nu}}-1)< r^{(i)}_{j,\nu_i} \le 2^{-j}k_{i,\widetilde{\nu}}$. The atomic decomposition theorem and \eqref{www1} give us
		\begin{equation}
		\label{op-below3}
		\big\|\psi_{j,\widetilde{\nu}}|H_2^{s/2}(\R^d)\big\|\;\le\; C_1  \left(\prod_{i=1}^m k_{i,\widetilde{\nu}}^{\gamma_i -1}\right)^{1/2} . 
		\end{equation}
		
		On the other hand  direct calculations show that the measure of the set $A_{j,\widetilde{\nu}}$ is equivalent to  $2^{-jd}\prod_{i=1}^m k_{i,\widetilde{\nu}}^{\gamma_i -1}$. So 
		\begin{eqnarray}\label{op-below4}		
		\big\|\psi_{j,\widetilde{\nu}}\sqrt{V}|L_2({\R}^d)\big\|&>& 2^{-js/2} 2^{j\frac d 2}\Big(\int_{A_{j,\widetilde{\nu}}}V(x)dx\Big)^{1/2}\;\ge\\
		&&C\, \varepsilon  2^{-js/2} 2^{j\frac d 2} | A_{j,\widetilde{\nu}}|^{1/2}\;\ge\;
		C_2 2^{-j\frac{s}{2}} \left(\prod_{i=1}^m k_{i,\widetilde{\nu}}^{\gamma_i -1}\right)^{1/2}
		.\nonumber
		\end{eqnarray}
		
		{\em Substep 2.2.} We choose $j=\big[ s^{-1}\log_2(C_1^{-1}C_2^2\beta)\big]$. Inequalities  \eqref{op-below3} and \reff{op-below4} imply
		\begin{align}\nonumber
		\big( \Delta^{s/2}_\theta \psi_{j,\widetilde{\nu}}\, , \, \psi_{j,\widetilde{\nu}}\big)& \le  C_1\left(\prod_{i=1}^m k_{i,\widetilde{\nu}}^{\gamma_i -1}\right) \;\le\;
		C_1C_2^{-2}2^{sj}\big\|\psi_{j,\widetilde{\nu}}\sqrt{V}|L_2({\R}^d)\big\|^2\;<\;\\
		& \beta \big(V \psi_{j,\widetilde{\nu}}\, , \, \psi_{j,\widetilde{\nu}}\big). \label{op-below5}\nonumber 
		\end{align}
		
		The subspace  $M={\rm span}\{\,\psi_{j,\widetilde{\nu}}\}_{\widetilde{\nu}}$ has the dimension  $\dim M= 2^{jm}\,\sim \, \beta^{m/s}$. The functions $\psi_{j,\widetilde{\nu}}$ are pairwise orthogonal therefore  for any $\psi\in M$ we have
		$$
		\big( H_{s,\beta,\theta} \psi\, , \, \psi\big)\;< \; 0\, .
		$$
		For any subspace $N\subset R_\gamma L_2({\R}^d)$ of dimension $\dim M -1$ one can find a function $\psi\in M$  such that $\|\psi\|=1$ and $\psi\perp N$. In consequence
		$$
		\sup_N \inf_{\psi\in {\mathcal D}(H_{s,\beta,\theta}), \|\psi\|=1, \psi\perp N} \big( H_{s,\beta,\theta} \psi\, , \, \psi\big)\;<\; 0\, ,
		$$
		where the supremum is taken over all $ M -1$ dimensional subspaces of $R_\gamma L_2({\R}^d)$.
		So the Max-Min principle implies that $H_{s,\beta,\theta}$ has at least $\dim M\sim \beta^{m/s}$
		negative eigenvalues, cf. eg. \cite[p.489]{EE}.  
		\epr 
		
		\begin{Rem}{\rm 
				1) Since $s>\frac{d}{r}$  and $\frac{m}{s}< \frac{rm}{d}\le \frac{r}{\gamma_1} $ we have always the gap between the upper and lower estimates. The  estimate is more precise  for small values of $s$ .  
				
				2) Apart of the radial case we have no block-radial functions satisfying the assumption of  Theorem \ref{mainT1} and Theorem \ref{op-main} in dimensions   $d=2$ and $d=3$. So the smallest possible dimension for which the assumption of Theorem \ref{op-main} are satisfied  is  $d=4$ with $\gamma_1=\gamma_2=2$. 
	
				 3) If $d\ge 4$ and $V\in L^r(\R^d)$ is the radial potential satisfying the assumption of  Theorem \ref{op-main} then the operator 
				$H_{s,\beta,\theta}$ has asymptotically at most $\beta^{r/d}$ eigenvalues with radial eigenfunctions. On the other hand such potential is $SO(\gamma)$-invariant for any $\gamma$.  So choosing $\gamma$ such that   $\min_i\gamma_i\ge 2$  we have  asymptotically at least  $\beta^{m/s}$ eigenvalues with $SO(\gamma)$-invariant eigenfunctions. If $\frac d r < s < [d/2]\frac d r$ then the operator      $H_{s,\beta,\theta}$ has eigenfunctions that are block-radial but not radial.   
				}
		\end{Rem}


		

		
		\appendix
		\section{Atomic decomposition for subspaces of invariant functions} 
		
		We recall the main idea of the method of atomic decomposition, and we refer the reader to  \cite{LS} where the atomic decompositions  related to the action of compact groups of isometries   are described in  details.   We assume that $\gamma_i\ge 2$ for any $i=1,\ldots, m$.
		
		We start with the following notions of  separations  and discretizations. 
		\begin{Def}\label{sep}
			Let $\varepsilon>0$ be a positive number, $\alpha=1,2,\ldots$ be a 
			positive integer  and $X$ a nonempty subset of $\R^d$. \\
			(a) A subset ${\mathcal H}$ of $X$ is said to be an $\varepsilon$--separation of $X$, 
			if the distance between any two distinct points of 
			$\mathcal H$ is greater than or equal to $\varepsilon$.  \\ 
			(b) A subset $\mathcal H$ of $X$ is called an  $(\varepsilon,\delta)$--discretization 
			of $X$ if it is an $\varepsilon$--separation of $X$ and   
			\[
			X\subset \bigcup_{x\in{\mathcal H}}B(x,\delta\varepsilon).
			\]
		\end{Def}

		\begin{Rem}\label{wwwr}{\rm 
				Let $n$ be a positive integer. If $\mathcal H$ is an 
				$(\varepsilon,\delta)$--discretization of $\R^d$ and $n\geq \delta$, then 
				the family $\{B(x,n\varepsilon)\}_{x\in {\mathcal H}}$ is an uniformly 
				locally finite covering of $\R^d$ with multiplicity that can be estimated from 
				above by a constant depending on $d$ and $n$, but independent of $\varepsilon$.}
		\end{Rem}
		
		Now we describe  discretizations related to the group $SO(\gamma)$. In this case we can proceed in the following way.  Let 
		$\{ x^{(j,i)}_{k,\ell}\}$, $k\in \N_0$, and $\ell=0,\ldots ,k^{\gamma_i-1}$ be 
		a $(2^{-j},\delta_i)$--discretization in $\R^{\gamma_i}$ related to the action of special orthogonal group $SO(\gamma_i)$ on  $\R^{\gamma_i}$.   We refer to \cite{SS1} for the construction of this type of discretitation. In particular we have 
		\begin{align*}
		\big|x^{(j,i)}_{k,\ell}\big| \sim k\, 2^{-j}  \quad \text{and}\quad x^{(j,i)}_{k,\ell}\in SO(\gamma_i)\cdot  x^{(j,i)}_{k,0}, \quad \ell= 0,\ldots k^{\gamma_i-1} .
		\end{align*}
		
		We put
		\begin{align*}
		{\mathcal H}_j& = \{ x^{(j)}_{\tilde{k},\tilde{l}}=(x^{(j,1)}_{k_1,l_1},\ldots ,
		x^{(j,m)}_{k_m,l_m}):\quad \tilde{k}=(k_1,\ldots k_m),\quad \tilde{\ell} = (k_1,\ldots k_m) 
		\},\\
		\intertext{and}
		\delta& = \sqrt{d}\max_i\{\delta_i\}. 
		\end{align*}
		The set ${\mathcal H}_j$ is a $(2^{-j},\delta)$--discretization of $\R^d$. Please note that   $\tilde{\ell}=0$ if $\tilde{k}=0$ and that $x^{(j)}_{0,0}$ is the origin.  
		If $\tilde{k}\neq 0$  then 
		\begin{equation}
		SO(\gamma)\big(x^{(j)}_{\tilde{k},\tilde{l}}\big) = \prod_{i:k_i\not= 0} SO(\gamma_i)\big(x^{(j)}_{k_i,l_i}\big) 
		\end{equation}
		and 
		\begin{equation}\label{www1}
		{\rm card}\Big( {\mathcal H}_j \cap SO(\gamma)\big(x^{(j)}_{\tilde{k},\tilde{l}}\big)\Big) = \prod_{i:k_i\not= 0} k_i^{\gamma_i-1}= C(j,\tilde{k}). 
		\end{equation}

		We assume that $s>0$ and $1\le p\le \infty$. The function $a_{j,\tilde{k},\tilde{\ell}}$ is called an $(s,p)$-atom centred at the point $x_{j,\tilde{k},\tilde{\ell}}\in {\mathcal H}_j$ if: 
		\begin{align} \
		\supp a_{j,\tilde{k},\tilde{\ell}} \, &\subset \, B(x^{(j)}_{\tilde{k},\tilde{\ell}}, 2\delta 2^{-j})\, , 
		\nonumber\\ 
		\sup_{y\in \R^d}|\partial^\alpha a_{j,\tilde{k},\tilde{\ell}}(y)| &\leq 2^{-j(s-|\alpha|-\frac{d}{p})} , \qquad |\alpha|\leq s+1 .\nonumber 
		\end{align}
		
		Let $f \in R_\gamma B^{s}_{p,q}(\R^d)$.  The atomic decomposition theorem asserts that any function $f \in B^{s}_{p,q}(\R^d)$ can be decomposed in the following way 
		\begin{align}\label{atomic0}
		&f = 
		\sum_{j=0}^\infty \, \sum_{\tilde{k}\in \N_0^m} 
		\sum_{\tilde{\ell}} s_{j,\tilde{k}}\, a_{j,\tilde{k},\tilde{\ell}},  \qquad (\text{convergence in }\quad {\mathcal S}')
		\intertext{with} 
		&\bigg(\sum_{j=0}^{\infty} \, \Big(\sum_{\tilde{k}\in \N_0^m} \, C(j,\tilde{k})\, |s_{j,\tilde{k}},|^{p} \Big)^{q/p}\bigg)^{1/q} 
		< \infty \, , \qquad s_{j,\tilde{k}}\in \C \label{atomic00}
		\end{align}
		(usual change if $q=\infty$). On the other hand any distribution represented by \eqref{atomic0} with \eqref{atomic00} belongs to $R_\gamma B^{s}_{p,q}(\R^d)$. Moreover the infimum over all possible representations of the expressions \eqref{atomic00} give us an equivalent norm in $R_\gamma B^{s}_{p,q}(\R^d)$.
		For the proof we refer to \cite[Theorem 2]{LS}, cf. also Remark 7 ibidem.

		\section{Wavelet characterizations of  Besov spaces with ${\mathcal A}_\infty$
			weights}
		A locally integrable function $w:\R^d \rightarrow \R_+$ belongs to the class ${\mathcal A}_\rho$,  $1<\rho<\infty$ if it satisfies the inequality 
		$$
		\frac{1}{|Q|}\int_Q w(x) dx \; \left(\frac{1}{|Q|}\int_Q w(x)^{-\rho'/\rho} dx \right)^{\rho/\rho'} \le A<\infty
		$$
		for all cubes $Q$ in $\R^d$. The class $\mathcal{A}_\infty$ is the union of all the ${\mathcal A}_\rho$ classes. 
		We recall briefly the wavelet characterization of the weighted Besov spaces  proved in \cite{HS}. Further information and  references concerning the wavelet theory can be found there.

		Let $\widetilde{\phi}$ be an orthogonal scaling function 
		on $\R$ with compact support and of sufficiently high regularity.
		Let $\widetilde{\psi}$ be an associated wavelet. Then the 
		tensor-product ansatz yields
		a scaling function $\phi$  and associated wavelets
		$\psi_1, \ldots, \psi_{2^{d}-1}$, all defined now on $\Rd$. 
		We suppose
		\[
		\widetilde{\phi} \in C^{N_1}(\R) \qquad \text{and} \qquad 
		\supp \widetilde{\phi} \subset [-N_2,\, N_2]
		\]
		for certain natural numbers $N_1$ and $N_2$.
		This implies
		\begin{equation}
		\label{2-1-2}
		\phi, \, \psi_i \in C^{N_1}(\R) \quad \text{and} \quad 
		\supp \phi ,\, \supp \psi_i \subset [-N_3,\, N_3]^d \, , 
		\quad i=1, \ldots \, , 2^{n}-1 \, .
		\end{equation}
		We shall use the standard abbreviations 
		\begin{equation}
		\label{convention}
		\phi_{\nu,k}(x) =  2^{\nu d/2} \, \phi(2^\nu x-k) \quad
		\text{and}\quad
		\psi_{i,\nu,k}(x) =  2^{\nu d/2} \, \psi_i(2^\nu x-k)\, . 
		\end{equation}
		
		\begin{T}\label{waveweight}  
			Let $0 < p,q \le \infty$ and let $s\in \R$. 
			Let $\phi$ be a scaling function and let $\psi_i$,  $i=1, \ldots ,2^d-1$, be
			the corresponding wavelets satisfying \eqref{2-1-2}. We assume that  $|s|<N_1$. 
			Then a distribution $f \in \mathcal{S}'(\Rd)$ belongs to $B^s_{p,q} (\Rd,w)$,  if, and only if,
			\begin{align*}\label{wavelet1}
			\| \, f \, |B^s_{p,q}(\Rd, w)\|^\star  = &  \Big\| \left\{\langle
			f,\phi_{0,k}\rangle \right\}_{k\in \Z^d} | \ell_p(w)\Big\| \\ 
			& + \sum_{i=1}^{2^d-1}
			\Big\| \left\{\langle f,\psi_{i,\nu,k}\rangle \right\}_{\nu\in \N_0, k\in \Z^n} | b^\sigma_{p,q}(w)\Big\|
			< \infty \, ,\nonumber
			\end{align*}
			where $\sigma=s+\frac d 2 - \frac d p$. 
			Furthermore,
			$\| \, f \, |B^s_{p,q}(\Rd, w) \|^\star $ may be used as an 
			equivalent $($quasi-$)$ norm in  $B^s_{p,q}(\Rd, w)$ and the map .
			\begin{equation} \nonumber
			\mathcal{W}: B^s_{p,q}(\Rd, w)\ni f\mapsto \Big(\langle f,\phi_{0,k}, \rangle, \langle f,\psi_{i,\nu,k}\rangle \Big)_{(i, \nu,k)}\in  l_p(w)\oplus b^\sigma_{p,q}(w)
			\end{equation}
			is a topological linear  isomorphism. 
		\end{T} 
		



Alicja Dota

Institute of Mathematics, Pozna\'n University of Technology, 

Pozna\'n, Poland

e-mail: alicja.dota@put.poznan.pl

\vspace{1cm}
Leszek Skrzypczak

Faculty of Mathematics and  Computer Science, 

Adam Mickiewicz University, Pozna\'n, Poland

e-mail: lskrzyp@amu.edu.pl

\end{document}